\documentclass[reqno]{amsart}
\usepackage{enumerate}
\usepackage{color}
\usepackage{mathrsfs}
\usepackage{hyperref}
\usepackage{float}
\usepackage{graphicx}%
\hypersetup{
colorlinks,
citecolor=blue,
filecolor=black,
linkcolor=blue,
urlcolor=magenta
}
\hypersetup{linktocpage}

\usepackage{scrextend}
\changefontsizes[16pt]{11pt}

\newtheorem{theorem}{Theorem}[section]
\newtheorem{lemma}[theorem]{Lemma}

\theoremstyle{definition}

\newtheorem{remark}[theorem]{Remark}

\newtheorem{assumption}[theorem]{Assumption}

\newcommand{\norm}[1]{\left\Vert#1\right\Vert}

\numberwithin{equation}{section}
\usepackage{amsfonts}

\usepackage{amsmath}
\usepackage{amssymb}

\begin{document}
\font\nho=cmr10
\def\dive{\mathrm{div}}
\def\cal{\mathcal}
\def\L{\cal L}

\def \ud{\underline }
\def\id{{\indent }}
\def\f{\frac}
\def\non{{\noindent}}
 \def\le{\leqslant} 
 \def\leq{\leqslant}
 \def\geq{\geqslant} 
\def\rar{\rightarrow}
\def\Rar{\Rightarrow}
\def\ti{\times}
\def\i{\mathbb I}
\def\j{\mathbb J}
\def\si{\sigma}
\def\Ga{\Gamma}
\def\ga{\gamma}
\def\ld{{\lambda}}
\def\Si{\Psi}
\def\f{\mathbf F}
\def\r{\hro{R}}
\def\e{\cal{E}}
\def\B{\cal B}
\def\A{\mathcal{A}}
\def\p{\mathbb P}

\def\tet{\theta}
\def\Tet{\Theta}
\def\hro{\mathbb}
\def\ho{\mathcal}
\def\P{\ho P}
\def\E{\mathcal{E}}
\def\n{\mathbb{N}}
\def\M{\mathbb{M}}
\def\dMu{\mathbf{U}}
\def\dMcs{\mathbf{C}}
\def\dMcu{\mathbf{C^u}}
\def\vk{\vskip 0.2cm}
\def\td{\Leftrightarrow}
\def\df{\frac}
\def\Wei{\mathrm{We}}
\def\Rey{\mathrm{Re}}
\def\s{\mathbb S}
\def\l{\mathcal{L}}
\def\C+{C_+([t_0,\infty))}
\def\o{\cal O}
\def\grad{\operatorname{grad}}

\begin{center}
{\LARGE\bf Generalized gravitational fields and well-posedness of the Boussinesq systems on non-compact Riemannian Manifolds}

{\bf Tran Thi Ngoc}\footnote{Faculty of Fundamental Sciences, East Asia University of Technology, Trinh Van Bo Street, Nam Tu Liem, Hanoi, Vietnam. Email: ngoctt@eaut.edu.vn}
and
{\bf Pham Truong Xuan}\footnote{Faculty of Mathematics and Informatics, Hanoi University of Science and Technology, 1 Dai Co Viet, Hanoi, Vietnam. Email: phamtruongxuan.k5@gmail.com}
\end{center}

{\bf Abstract.} We study the global existence, uniqueness and exponential stability of mild solutions to the Boussinesq systems equipped with a generalized gravitational field on the framework of non-compact Riemannian manifolds. We work on some manifolds satisfying some bounded and negative conditions on curvature tensors. We consider a couple of Stokes  and heat semigroups associated with the corresponding linear system which provides a vectorial matrix semigoup. By using dispersive and smoothing estimates of the vectorial matrix semigroup we establish the global-in-time existence and uniqueness of mild solutions for linear systems. Next, we can pass from the linear system to the semilinear systems to obtain the well-posedness by using fixed point arguments. Moreover, we will prove the exponential stability of such solutions by using Gronwall's inequality.   




\tableofcontents

\section{Introduction}
In the present paper, we are concerned with the incompressible Boussinesq system on non-compact Riemannian manifolds $(\mathbf{M},g)$ which satisfies $(H_1)-(H_4)$ conditions in Assumption \ref{Ass} below (where we consider that the dimension of manifold is $d\geqslant 2$ and $g$ is the Riemannian metric):
\begin{equation}
\left\{
\begin{array}
[c]{rll}%
u_{t}+(u\cdot\nabla)u - Lu+\nabla p & = \kappa\theta h + F\quad &
x\in\mathbf M,\,t>0,\hfill\\
\operatorname{div}u\!\! & =\;0\quad & x\in\mathbf M,\,t\geq0,\\
\theta_{t}-\widetilde{L}\theta+(u\cdot\nabla)\theta\!\! & = f\quad & x\in
\mathbf M,\,t>0,\\
u(x,0)\!\! & =\;u_{0}(x)\quad & x\in\mathbf M,\\
\theta(x,0)\!\! & =\;\theta_{0}(x)\quad & x\in\mathbf M,
\end{array}
\right.  \label{BouE}%
\end{equation}
where $L= \overrightarrow{\Delta}+r$ is Ebin-Marsden's Laplace operator defined by the stress tensor (where $r$ is the Ricci operator), $\widetilde{L}= \Delta_g$ is Laplace-Beltrami operator associated with metric $g$, the constant $\kappa>0$ is the volume expansion coefficient. The field $h(t,x)$ is a generalization of the gravitational field on manifold ${\bf M}$ and the constant $\kappa>0$ is the volume expansion coefficient. The unknowns: $u(t,x)$ is the velocity field, $p(t,x)$ is the scalar pressure, and the scalar function $\theta(t,x)$ is the temperature. The vector fields $f$ and $F$ represent the reference temperature and the external force, respectively (for details see Section \ref{S2}). 
Considering the zero-temperature case, i.e., $\theta=0$, then system
(\ref{BouE}) becomes the Navier-Stokes equations.

We now reall briefly some results on the Boussinesq system in Euclidean space $\mathbb{R}^d$. Fife and Joseph \cite{Fi1969} provided one of the first rigorous mathematical results for the
convection problem by constructing analytic stationary solutions for the Boussinesq system with the bounded field $h$, as well as analyzing some stability and
bifurcation properties. After, Cannon and DiBenedetto \cite{Ca1980}
established the local-in-time existence in the class $L^{p}(0,T;L^{q}%
(\mathbb{R}^{n}))$ with suitable $p,q$. Hishida \cite{Hi1995} (see also \cite{Mo1991}) obtained the existence and exponential stability of global-in-time strong solutions for the Boussinesq system near to the steady state in a bounded domain of $\mathbb{R}^{3}$. Later, by using the $L^{p,\infty}$-$L^{q,\infty}$-dispersive and smoothing estimates in weak-$L^{p}$ spaces of the semigroup $e^{-tL}$ associated with the corresponding linear equations of the Boussinesq system,
Hishida \cite{Hi1997} showed the existence and large-time behavior of
global-in-time strong solutions in an exterior domain of $\mathbb{R}^{3}$
under smallness assumptions on the initial data $(u_{0},\theta_{0})$.
Well-posedness of time-periodic and almost periodic small solutions in exterior domains were proved
in \cite{HuyXuan22, Na2020} by employing frameworks based on weak-$L^{p}$ spaces. The
existence and stability of global small mild solutions for the Boussinesq system were studied in weak-$L^{p}$ spaces in \cite{Fe2006,Fe2010} and in Morrey spaces in \cite{Al2011}. A  result of stability in $B_{2,1}^{3/2}\times \dot{B}_{2,1}^{-1/2},$ under small perturbations, for a class of global large
$H^{1}$- solutions was proved by \cite{Liu2014}. Brandolese and Schonbek \cite{Br2012} obtained results on the existence and time-decay of weak solutions for the Boussinesq system in whole space $\mathbb{R}^{3}$ with initial data $(u_{0},\theta_{0})\in L^{2}\times L^{2}$. Li and Wang \cite{Li-Wang2021}
analyzed the Boussinesq system in the torus $\mathbb{T}^{3}$ and obtained an
ill-posedness result in $\dot{B}_{\infty,\infty}^{-1}\times\dot{B}%
_{\infty,\infty}^{-1}$ by showing the so-called norm inflation phenomena. Komo
\cite{Komo2015} analyzed the Boussinesq system in general smooth domains
$\Omega\subset$ $\mathbb{R}^{3}$ and obtained uniqueness criteria for strong
solutions in the framework of Lebesgue time-spatial mixed spaces
$L^{p}(0,T;L^{q}(\Omega))$ by assuming $(u_{0},\theta_{0})\in L^{2}\times
L^{2}$ and $g\in L^{8/3}(0,T;L^{4}(\Omega))$. Considering the case of a
constant field $h$, Brandolese and He \cite{Br2020} showed the uniqueness of
mild solutions in the class $(u,\theta)\in C([0,T],L^{3}(\mathbb{R}^{3})\times
L^{1}(\mathbb{R}^{3}))$ with $\theta\in L_{loc}^{\infty}((0,T);L^{q,\infty
}(\mathbb{R}^{3}))$. The existence and uniqueness results in the partial inviscid cases of the Boussinesq system were studied in \cite{Danchin2009,Danchin2008}, where the authors
explored different kinds of conditions on the initial data $(u_{0},\theta
_{0})$ involving $L^{p},$ $L^{p,\infty}$ (weak-$L^{p}$) and Besov spaces.
Recently, the unconditional uniqueness of mild solutions for Boussinesq equations in Morrey-Lorentz spaces has established by Ferreira and Xuan \cite{FePha2023}. The existence and stability of periodic mild solutions of Boussinesq systems in weak-Morrey spaces have proven by Xuan et al. \cite{XVT2024}. In the same framework, the asymptotic behaviour of solutions of stationary Boussinesq systems have been established by Xuan and Ngoc in a recent paper \cite{XN2024}.

Concerning the fluid dynamic equations on Riemannian manifolds, we summarize some previous results about the well-posedness, ill-posedness, asymptotic behaviour etc... on compact and noncompact manifolds.
The study of Navier-Stokes equations associated with the vectorial laplace operator given by the stress  tensor formula (known as Ebin and Marsden's laplace operator) on noncompact Einstein manifolds with negative  Ricci curvature was mentioned initially by Ebin and Marsden \cite{EbiMa} when they formulated  the formula for Navier-Stokes equations on an Einstein manifold with negative Ricci curvature. Since then, this notion has been used in the works of Czubak and Chan \cite{Cz1,Cz2} and also Khesin and Visiolek \cite{Khe} to prove the non-uniqueness of weak Leray solution of Navier-Stokes equation on the two-dimensional  hyperbolic manifolds. In additionally, the non-uniqueness problem of Navier-Stokes equations were extended to study on another three-dimensional non-compact Riemannian manifold by Lichtenfelz \cite{Li2016} by considering Anderson's manifolds.
Furthermore, Pierfelice \cite{Pi} has proved the dispersive and smoothing estimates for Stokes semigroups on the generalized non-compact manifolds with negative Ricci curvature then combines these estimates with Kato-iteration method to prove the existence and uniqueness of strong mild solutions to Navier-Stokes equations. The existence and stability of periodic and asymptotically almost periodic mild solutions to the Navier-Stokes equations on noncompact manifolds with negative curvature tensors have been established in some recent works \cite{HuyXuan2021,HuyXuan2022,XVQ2023, XuanVan2021}. In the related works, the Navier-Stokes equations associated with Hodge-Laplace operator has been studied in several  manifolds, e.g., on two sphere \cite{Cao1999,Il1991}, on compact Riemannian manifolds \cite{Fa2018,Fa2020,Ko2008,MiTa2001,Sa}, or on the connected sums of $\mathbb{R}^3$ in \cite{Zha}. 

Notably, in recent works \cite{XuanTrung2022,NgocXuan2024}, Xuan et al. have studied the well-posedness and exponential stability of periodic and pseudo almost periodic mild solutions for Boussniesq systems on the real hyperbolic manifolds. In these works, the authors first provided a model of Boussinesq system equipped with a generalized gravitational field on the real hyperbolic manifold, i.e., a non-Euclidean framework. Then, they used the $L^p-L^q$-dispersive and smoothing estimates of the scalar heat and Stokes semigroups to establish the well-posedness of mild solutions for linear and semilinear systems. Also, they proved the Messara-type principles for Boussinesq systems which guarantee the existence of periodic and pseudo almost periodic mild solutions. The works \cite{XuanTrung2022,NgocXuan2024} can be extended to the non-compact Riemannian manifold $(\mathbf M,g)$ which satisfying the conditions in Assumption \ref{Ass} and this paper is the first step for establishing the well-posedness of Boussinesq systems in this framework.

We describle the strategy of this paper as follows: first, we represent system \eqref{BouE} under the matrix intergral equation (see equation \eqref{MatrixEq} below). Then, we use the estimates for the Stokes and scalar heat semigroups (obtained in \cite{Pi}) to prove the $L^p-L^q$-dispersive and smoothing estimates for the vectorial matrix semigroup asscociated with the corresponding linear system of Boussinesq system \eqref{BouEq} (see Lemma \ref{dispersive}). Using these estimates and the boundedness of gamma functions we prove the existence of bounded mild solution for the linear system (see Theorem \ref{thm1}).
Note that, the non-commutation of the Kodaira-Hodge operator and the Stokes semigroup leads us to consider the existence on the phase spaces with time weights and gradient functions as well as in \cite{Pi,HuyXuan2022}. This is also an extension from the works on hyperbolic spaces in \cite{XuanTrung2022}.
After that, we establish the estimates for the bilinear operator associated with Boussinesq system, i.e., bilinear estimates \eqref{Bilinear-2} and \eqref{Bilinear-p^}. Combining these estimates with the existence results of the linear system and fixed point arguments, we establish the existence of bounded mild solution for the Boussinesq system in Theorem \ref{PeriodicThm}. We use Gronwall's inequality and the boundedness of beta and gamma functions to prove the exponential stability of the Boussinesq system (see Theorem \ref{stability}). This paper together with \cite{XuanTrung2022, NgocXuan2024} are continuation of the study on fluid dynamic systems on non-compact Riemannian manifolds satisfying bounded geometric conditions in Assumption \ref{Ass}.

This paper is organized as follows: in Section \ref{S2}, we present non-compact Riemannian manifolds satisfying some bounded conditions on curvatures and the setting of Boussinesq systems on these manifolds; in Section \ref{S3}, we first provide the $L^p-L^q$-dispersive and smoothing estimates and the proofs of the well-posedness of linear and Boussinesq systems, then we prove the exponential stability of the mild solutions for the Boussinesq system; in Subsection \ref{appen}, we discuss more details about the gravitational fields on non-compact Riemannian manifolds and in Subsection \ref{convergence} we give some detailed calculations for the integrals's convergences.\\
{\bf Notations.} In order to the convenience, through this paper we use the following notations
\begin{itemize}
\item We denote $(L^p\cap L^q)(X):= L^p(X)\cap L^q(X)$ and  $\norm{\cdot}_p:=\norm{\cdot}_{L^p(X)}$ on the space $L^p(X)$;
\item We define and use the norm $\norm{\begin{bmatrix}
u\\ \theta
\end{bmatrix}}_{L^p\times L^p} = \max\left\{ \norm{u}_{L^p({\mathbf M};\Gamma(T{\mathbf M}))},\norm{\theta}_{L^p({\mathbf M};\mathbb{R})}  \right\}$ on the Cartesian product space $ L^p({\mathbf M};\Gamma(T{\mathbf M})) \times  L^p({\mathbf M};\mathbb{R}) $;

\item ${\bf B}(\cdot,\cdot)$ and ${\bf\Gamma}(\cdot)$ denote beta and gamma functions respectively.
\end{itemize}  
{\bf Acknowledgments:} P. T. Xuan was funded by the Postdoctoral Scholarship Programme of Vingroup Innovation Foundation (VINIF), code VINIF.2023.STS.55

\section{Boussinesq system on non-compact manifolds}\label{S2}
\subsection{Non-compact Riemannian manifolds}\label{NSERiemann}
\noindent
 
In this subsection, we recall some notions on differential operators on Riemannian manifolds. We refer the reader to \cite{He,Jo} for notions and  detail discussions on Riemannian manifolds and related concepts of geometric analysis.

Let $(\mathbf M,g)$ be a $d$-dimensional Riemannian manifold (where $d\geqslant 2$), we denote the Levi-Civita connection by $\nabla$ and the set of all vector fields on $\mathbf M$ by $\Gamma(T\mathbf M)$. For $X\in \Gamma(T\mathbf M)$ we can extend  $\nabla_X$ to arbitrary $(p,q)$
tensor by requiring
\begin{enumerate}
	\item $\nabla_X(c(S)) = c(\nabla_X S)$ for any contraction $c$,
	\item $\nabla_X(S \otimes T) = \nabla_X S \otimes T + S \otimes \nabla_X T$
\end{enumerate}   
where we take the convention that $\nabla_X f = X\cdot f$ for a function $f:\mathbf M\to \r$. 

In particular, for 
$S\in \Gamma(\otimes^p(T\mathbf M)\otimes^q
(T^*\mathbf M))$  we get 
$$(\nabla_X S)(X_1, \cdots, X_q ) = \nabla_X(S(X_1, \cdots, X_q ))-
S(\nabla_XX_1,\cdots, X_q))-\cdots-S(X_1, \cdots,\nabla_X X_q).$$
Moreover, we define the covariant derivatives $\nabla$ on tensor field 
$S\in \Gamma(\otimes^p(T\mathbf M)\otimes^q(T^*\mathbf M))$ 
by
$$\nabla S(X, X_1, \cdots, X_q) = (\nabla _X S)(X_1, \cdots, X_q),
\hbox{ hence, }\nabla S\in \Gamma(\otimes^p(T\mathbf M)\otimes^{q+1}(T^*\mathbf M).$$
Next, we recall the "music notations" on Riemannian manifolds. 
For a 1-form $w$, we define the vector field $w^\sharp$	 by
$$g(w^\sharp, Y ) = w(Y), \forall Y \in \Gamma(T\mathbf M)$$
whereas, for a vector field $X$, we define the 1-form $X^\flat$ by
$$X^\flat(Y) = g(X, Y ), \forall Y \in \Gamma(T\mathbf M).$$
The metric on $1$-forms  can then be defined by setting
$$g(w, \eta) := g(w^\sharp, \eta^\sharp), \forall w, \eta \in \Gamma(T^*\mathbf M).$$
The Riemannian gradient of a function is then defined as
$$\grad p = (dp)^\sharp.$$
More generally, for $(p,q)$-tensor field $ T\in \Gamma(\otimes^p(T\mathbf M)\otimes^q(T^*\mathbf M))$
we have
\begin{eqnarray*}
	T^\sharp& = &C^2_1(g^{-1}\otimes T) \in \Gamma(\otimes^{p+1}(T\mathbf M)\otimes^{q-1}(T^*\mathbf M)),\cr
	T^\flat & = &C^1_2(g \otimes  T ) \in \Gamma(\otimes^{p-1}(T\mathbf M)\otimes^{q+1}(T^*\mathbf M)),\cr
	\dive T &= &C^1_1\nabla T \in \Gamma(\otimes^{p-1}(T\mathbf M)\otimes^{q}(T^*\mathbf M))
\end{eqnarray*}
where $C^i_j$ stands for the contraction of the $i$ and $j$ indices for tensors.

The Riemann curvature tensor $\mathcal{R}: \Gamma(T{\bf M})\times \Gamma(T{\bf M}) \times \Gamma(T{\bf M}) \to \Gamma(T{\bf M})$ is defined by 
$$\mathcal{R}(X, Y )Z :=
-\nabla_X (\nabla_Y Z)+
\nabla_Y (\nabla_X Z)+\nabla_{[X,Y ]}Z \hbox{ for all }X, Y, Z  \in \Gamma(T\mathbf M).$$
Then, the Ricci curvature tensor $\mathrm{Ric}: \Gamma(T{\bf M})\times \Gamma(T{\bf M}) \to \Gamma(T{\bf M})$ is defined from $\mathcal{R}$ as
$$\mathrm{Ric}(X,Y)=\sum_{i=1}^dg(\mathcal{R}(X,e_i)Y,e_i)\hbox{ for all }X, Y\in\Gamma({T\mathbf M})$$
for the standard basis $\left\{e_i=\dfrac{\partial}{\partial x_i}\right\}_{i=1}^d$. 
Moreover, the sectional
curvature $\kappa$ is given by 
$$\kappa(X,Y):=\frac{\cal R(X,Y,X,Y)}{g(X,X)g(Y,Y)-g(X,Y)^2} $$
for all vector fields $X, Y\in T_x\mathbf M$.

In the rest of this paper, we work on a  smooth, complete, non-compact, simply connected Riemannian manifold $(\mathcal{M},g)$ which satisfies the following hypotheses (which were considered  by Pierfelice in \cite{Pi}):
\begin{assumption}\label{Ass}
\noindent

\begin{itemize}
	\item[$(H_1)$] $|\mathcal R| + |\nabla \mathcal R| + |\nabla^2 \mathcal R| \leq K \, ,$
	\item[$(H_2)$] $-\dfrac{1}{c_0}g \leq \mathrm{Ric} \leq -c_0 g \, ,$ for some $c_0$ positive,
	\item[$(H_3)$] $\kappa <0$,
	\item[$(H_4)$] $\inf_{x\in \mathcal{M}} r_x >0$,
\end{itemize}
where 
$r_x$ is the injectivity radius for the exponential map at $x$.
\end{assumption}

There are several Riemannian manifolds satisfying the above hypotheses $(H_1)-(H_4)$ such as 
real hyperbolic manifolds,  non-compact Einstein manifolds with negative Ricci curvature tensors  (\cite{Hel,Jo}), Damek-Ricci manifolds (\cite{Da}) and symmetric manifolds of non-compact types (\cite{Er, Hel}).

We also recall the Laplace-Beltrami operator $\Delta_g$ applying on functions which is defined as 
	$$\Delta_g(f)=\dive\grad f =\frac{1}{\sqrt{|g|}}\frac{\partial}{\partial x^j}\left(\sqrt{|g|}g^{ij}\frac{\partial f}{\partial x^i}\right)\hbox{ for a function }f :\mathbf {M} \to \r,$$
	where $|g| = \mathrm{det}\,g$.

Furthermore,  the vectorial Laplacian $L$ is defined by the stress tensor (see \cite{Pi} and also \cite{EbiMa,Tay}):
$$Lu = \dive(\nabla u + \nabla u^t)^{\sharp}.$$
Since $\dive u=0$ we can express $L$ in the following formula
$$Lu = \overrightarrow{\Delta}u + r(u),$$
where $\overrightarrow{\Delta}$ is the Bochner-Laplacian
$$\overrightarrow{\Delta}u =- \nabla^*\nabla u= \mathrm{Tr}_g(\nabla^2u)$$
and $r(\cdot)$ is the Ricci operator related to the Ricci curvature tensor by the formula
$$r(u) = (\mathrm{Ric}(u,\cdot))^{\sharp} \hbox{ for all }u\in\Gamma({T\mathbf M})$$.

\subsection{Boussinesq systems and generalized gravitational fields}
For the sake of simplicity, we will deal with the incompressible Boussinesq system on non-compact manifolds ${\bf M}$ with the volume expansion coefficient $\kappa=1$:
\begin{equation}\label{BouEq} 
\left\{
  \begin{array}{rll}
 u_t + (u\cdot\nabla)u - L u + \nabla p \!\! &= \theta h + F, \hfill \\
\nabla \cdot u \!\!&=\; 0, \\
\theta_t - \widetilde{L} \theta + (u \cdot \nabla)\theta \!\!&= f, \\
u(0) \!\!& = \;u_0,\\
\theta(0) \!\!& = \;\theta_0,\\
\end{array}\right.
\end{equation}
where, the unknowns are $u(x,t): {\bf M} \times \r \to \Gamma(T{\bf M}), p(x,t): {\bf M} \times \r \to \mathbb{R}$ and $\theta(x,t): {\bf M} \times \r \to \mathbb{R}$ representing  the velocity field, the pressure and the temperature of the fluid at point $(x,t) \in  {\bf M}\times \r$, respectively.
Here, the vectorial Laplce operator $L= \overrightarrow{\Delta}+r$ is defined in Subsection \ref{NSERiemann}, the scalar operator $\widetilde{L}= \Delta_g$ is Laplace-Beltrami operator associated with metric $g$.
The functions $f: {\bf M}\times \mathbb{R} \to \mathbb R$ is provided such that $ f$ represents the reference temperature and $F: {\bf M} \times \mathbb{R} \to \Gamma(T{\bf M})$ is a second order tensor fields such that $F$ represents the external force. The funtion $h(x):{\bf M}\to {\Gamma}(T{\bf M})$ represents the gravitational field

Normally, the funtion $h$ (depends only on variable $x\in {\bf M}$) represents the gravitational field and does not depend on the time. However, in this work we are going to study the general case, which $h(x,t):{\bf M}\times \mathbb{R} \to \Gamma(T{\bf M})$ depends on the time line and satisfies the following assumption which guarantees the regularity for elliptic problem to determine the pressure $p$:
\begin{assumption}\label{Assum}
Assume that function $h(\cdot,t)$ satisfies
\begin{equation}
h\in C_b(\mathbb{R}_+, L^\infty(\Gamma(T{\bf M}))) \hbox{   and   } h\in C_b(\mathbb{R}_+, L^{\frac{d}{2},\infty}(\Gamma(T{\bf M}))).
\end{equation}
\end{assumption}
The existence of generalized gravitational fields satisfied Assumption \ref{Assum} will be given in Subsection \ref{appen}.

Now, we need to fixed a pressure $p$ to guarantees the uniqueness of Leray weak solutions of Boussinesq system \eqref{BouEq}. By the same way as in \cite{Pi} (see Section 3), we take divergence to the first equation of system \eqref{BouEq} and using $\dive(\overrightarrow{\Delta}u) =r(u)$ (if $\dive u=0$), we get
\begin{equation}\label{EllipticEq}
\Delta_g p = \dive[2r(u)-(u\cdot \nabla)u + \theta h + F]. 
\end{equation} 
If we consider $u \in C_b(\r_+,L^{p}({\bf M};\Gamma(T{\bf M})))$, $\theta \in C_b(\r_+,L^p({\bf M}; \mathbb{R}))$, $h \in C_b(\r_+, L^{\infty}({\bf M};\Gamma(T{\bf M})))$ and $F \in C_b(\r_+,L^{p}({\bf M};\Gamma(T{\bf M}\otimes T{\bf M})))$, then we have 
$$2r(u)-(u\cdot \nabla)u + \theta h + F  = 2r(u)-\dive(u\otimes u) + \theta h + F \in C_b(\r_+, L^p({\bf M};\Gamma(T{\bf M}))).$$ 
Here, we used $(u\cdot \nabla)u = \dive (u\otimes u)$ because $\dive u=0$ and $r(\cdot)$ is bounded by $(H_2)$ condition.

Moreover, by hypotheses $(H_1-H_4)$ we have that the operator $\Delta_g \,: \, W^{2,q}({\bf M};\mathbb{R}) \rightarrow L^q({\bf M};\mathbb{R})$ is an isomorphism for $2\leq q <\infty$ (see \cite{Loho}). Therefore, for $p>1$, we can choose the solution of elliptic equation \eqref{EllipticEq} by
\begin{equation}
p= \Delta_g^{-1}\dive[2r(u)-(u\cdot \nabla)u + \theta h + \dive F].
\end{equation}
This choice guarentees to get the uniqueness of Leray weak solutions for Boussinesq system (see also \cite[Section 6]{Pi} for the uniequeness of Navier-Stokes equations). From above formula, we have 
\begin{equation}\label{nap}
\nabla p= \nabla(-\Delta_g)^{-1}\dive[-2r(u)+ (u\cdot \nabla)u - \theta h - \dive F].
\end{equation}
Since Riesz transforms are $L^p$-bounded on ${\bf M}$ (see \cite{Loho}), we obtain that the operator $\nabla(-\Delta_g)^{-1} \mathrm{div}:L^p({\bf M};\Gamma(T{\bf M})) \to L^p({\bf M};\Gamma(T{\bf M}))$ is bounded. Therefore, we have $\nabla p\in L^p({\bf M};\Gamma(T{\bf M}))$.

By setting $\mathbb{P} = \nabla(-\Delta_g)^{-1} \mathrm{div} + \mathrm{Id}$ (which is called Kodaira-Hodge operator) and replacing the formula \eqref{nap} of $\nabla p$ into system \eqref{BouEq}, we recieve the following system
\begin{align}\label{AbstractE}
\begin{cases}
\partial_t u = (L + G)u + \mathbb{P}(\theta h) + \mathbb{P}[- (u\cdot\nabla)u + F],\\
\partial_t\theta = \widetilde{L}\theta -(u\cdot\nabla)\theta + f,\\
\dive{u}= 0, \\
u(0) = u_0,\, \theta(0) = \theta_0,
\end{cases}
\end{align}
where $Gu = -2\mathrm{grad}(-\Delta_g)^{-1}\mathrm{div}(ru)$. 

\begin{remark}\label{Rem1}
Note that, in general the operator $L+G= \overrightarrow{\Delta} + r + G$ is not commutative with the Kodaira-Hodge operator
	$\mathbb{P} = I + \mathrm{grad}(-\Delta_g)^{-1} \mathrm{div}$ on the non-compact manifolds which satisfy all the conditions $(H_1)-(H_4)$ (see \cite{Pi}). In a specific case, if the Ricci curvature tensor is a multiplication of the metric $g$ with a constant, i.e., Einstein manifolds, then $G=0$ and $L+G= \overrightarrow{\Delta} + r$ is commutative with $\mathbb{P}$.
\end{remark}

In what follows, one will investigate the system \eqref{AbstractE} with $(u,\theta)$ in the product space $C_b(\mathbb{R}_+,L^p({\bf M};\Gamma(T{\bf M})))\times C_b(\mathbb{R}_+,L^p({\bf M};\mathbb{R})).$ To do this, we need to put the new operator $\mathcal{A}:=%
\begin{bmatrix}
-(L+G) & 0\\
0 & -\widetilde{L}
\end{bmatrix}
$, which acts on the Cartesian product space $L^p({\bf M};\Gamma(T{\bf M}))\times L^p({\bf M};\mathbb{R})$. Hence,
by utilizing Duhamel's principle in a matrix form, we arrive at the following
integral formulation for (\ref{AbstractE})%
\begin{equation}\label{MatrixEq}
\begin{bmatrix}
u(t)\\
\theta(t)
\end{bmatrix}
=e^{-t\A}%
\begin{bmatrix}
u_{0}\\
\theta_{0}%
\end{bmatrix}
+\mathscr B\left(
\begin{bmatrix}
u\\
\theta
\end{bmatrix}
,%
\begin{bmatrix}
u\\
\theta
\end{bmatrix}
\right)  (t)+\cal H_h(\theta)(t) + \mathcal{F} \left(
\begin{bmatrix}
F\\
f
\end{bmatrix}
\right)(t), 
\end{equation}
where the bilinear, linear-coupling and external forced operators used in the above equation
are given respectively by
\begin{equation}
\mathscr B\left(
\begin{bmatrix}
u\\  \zeta
\end{bmatrix}
,%
\begin{bmatrix}
v\\  \vartheta
\end{bmatrix}
\right)  (t):=-\int_{0}^{t} e^{-(t-\tau)\A}
\begin{bmatrix}
\mathbb{P}[(u\cdot\nabla)v]\\
(u\cdot\nabla) \vartheta
\end{bmatrix}
(\tau)d\tau \label{Bilinear},
\end{equation}
 
\begin{equation}
\cal H_h(\vartheta)(t):=\int_{0}^{t}e^{-(t-\tau)\A}%
\begin{bmatrix}
\mathbb{P}[\vartheta h]\\
0
\end{bmatrix}
(\tau)d\tau \label{Couple}%
\end{equation}
and
\begin{equation}
\mathcal{F} \left(
\begin{bmatrix}
F\\
f
\end{bmatrix}
\right)(t):= \int_0^t e^{-(t-\tau)\A}
\begin{bmatrix}
\mathbb{P} [F]\\
f
\end{bmatrix}
(\tau)d\tau.   
\end{equation}

\section{The global-in-time well-posedness and exponential stability}\label{S3}
\subsection{Linear estimates and well-posedness for the linear systems}
 
One first studies to the following linear equation corresponding to the integral matrix equation \eqref{MatrixEq} as follows. 
\begin{equation}\label{LinearE}
\begin{bmatrix}
u(t)\\
\theta(t)
\end{bmatrix}
=e^{-t\A}%
\begin{bmatrix}
u_{0}\\
\theta_{0}%
\end{bmatrix}
+\cal H_h(\eta)(t) + \mathcal{F} \left( \begin{bmatrix}
F\\
f
\end{bmatrix} \right) (t).
\end{equation}
Here 
\begin{equation}
\cal H_h(\eta)(t):=\int_{0}^{t}e^{-(t-\tau)\A}%
\begin{bmatrix}
\mathbb{P}[\eta h] \\
0
\end{bmatrix}
(\tau)d\tau. \label{Couple}%
\end{equation}

Note that the $L^p-L^q$ dispersive and smoothing properties of the matrix semigroup $e^{-t\A}$  are important tools to estimate the bilinear, linear-coupling and external forced operators ($\mathscr B(\cdot,\cdot)$, $\cal H_h(\cdot)$ and $\mathcal{F}(\cdot)$, respectively)  in the equation \eqref{MatrixEq}. Therefore, we will show more clearly for this properties in the below lemma. 

\begin{lemma}\label{estimates} Suppose that Asumption \ref{Ass} holds. Then, there are $\beta \geq c_0>0$ and some $C>0$ such that the solution of \eqref{MatrixEq} satisfying
\item[$(i)$] For $t>0$, and $p$, $q$ such that $2\leq p \leq q < \infty$, the following dispersive estimates of vectorial matrix semigroup $e^{-t\mathcal{A}}$ hold
\begin{equation}\label{dispersive}
\left\| e^{-t \A} \begin{bmatrix}
u_0\\
\theta_0
\end{bmatrix}\right\|_{L^q} \leq C [h_d(t)]^{\frac{1}{p}-\frac{1}{q}}e^{-\beta t}\left\| \begin{bmatrix}
u_0\\
\theta_0
\end{bmatrix} \right\|_{L^p\cap L^2} 
\end{equation}
for all $(u_0,\theta_0) \in (L^p\cap L^2)({\mathbf M};\Gamma(T{\mathbf M})) \times (L^p\cap L^2)({\mathbf M};\mathbb{R})$, 
 where 
 $$h_d(t): = \max\left( \frac{1}{t^{d/2}},1 \right).$$   
\item[$(ii)$] For $>0$, and $p,q$ such that $2\leq p\leq q <\infty$ we have the following smoothing estimates of $e^{-t\mathcal{A}}$: 
\begin{equation}\label{smoothing}
\left\| \nabla e^{-t\A}\begin{bmatrix}
u_0\\ \theta_0
\end{bmatrix} \right\|_{L^q} \leq C [h_d(t)]^{\frac{1}{p}-\frac{1}{q}+\frac{1}{d}}e^{-\beta t} \left\| \begin{bmatrix}
u_0\\ \theta_0
\end{bmatrix} \right\|_{L^p\cap L^2}
\end{equation}
for all  $(u_0,\theta_0) \in (L^p\cap L^2)({\mathbf M};\Gamma(T{\mathbf M})) \times (L^p\cap L^2)({\mathbf M}; \mathbb R)$. 
\end{lemma}
\begin{proof}
$(i)$ First, we have
\begin{equation*}
 e^{-t\A} =  \begin{bmatrix}
e^{t(L+G)}&&0\\
0&& e^{t\widetilde{L}}
\end{bmatrix}.
\end{equation*}
The $L^p-L^q$ dispersive and smoothing estimates for the Stokes semigroup $e^{t(L+G)}$ and the scalar heat semigroup $e^{t\widetilde{L}}$ which were been establised by Pierfelice \cite{Pi}. In particular, for all $t>0$, and $p$, $q$ such that $2\leq p \leq q <\infty$, 
the $L^p-L^q$ dispersive estimates  of Stokes semigroup $e^{t(L+G)}$ is (see \cite[Corollary 4.13]{Pi} and its proof): 
\begin{equation}\label{di1}
\left\| e^{t (L+G)} u_0\right\|_{L^q} \leq C [h_d(t)]^{\frac{1}{p}-\frac{1}{q}}e^{-\beta t}\left\| u_0 \right\|_{L^p\cap L^2}
\end{equation}
for all $u_0 \in (L^p\cap L^2)({\mathbf M};\Gamma(T{\mathbf M}))$ and the $L^p-L^q$ dispersive estimates  of scalar heat semigroup $e^{t\widetilde{L}}$ is
\begin{equation}\label{di2}
\left\| e^{t\widetilde{L} } \theta_0\right\|_{L^q} \leq C [h_d(t)]^{\frac{1}{p}-\frac{1}{q}}e^{-\beta t}\left\| \theta_0 \right\|_{L^p\cap L^2} \hbox{ for all }\theta_0 \in (L^p\cap L^2)({\mathbf M};\mathbb{R}),    
\end{equation}
where $h_d(t): = \max\left( \frac{1}{t^{d/2}},1 \right)$. 
Clearly, the inequalities \eqref{di1} and \eqref{di2} lead to dispersive estimate \eqref{dispersive} for the semigroup $e^{-t\mathcal{A}}$.

$(ii)$  We have the gradient estimate \eqref{smoothing} in Assertion $(ii)$ since the fact that: for $t>0$ and $2\leq p\leq q <\infty$ the following $L^p-L^q$ smoothing estimates of Stokes semigroup $e^{t(L+G)}$ and of scalar heat semigroup $e^{t\widetilde{L}}$ hold (see \cite[Theorem 4.15]{Pi} and its proof):
\begin{eqnarray*}
\left\| \nabla e^{t(L+G)} u_0 \right\|_{L^q} &\leq&C [h_d(t)]^{\frac{1}{p}-\frac{1}{q}+\frac{1}{d}}e^{-\beta t} \left\| u_0\right\|_{L^p\cap L^2}
\end{eqnarray*}
and
\begin{equation*}
\left\|\nabla e^{t\widetilde{L}}\theta_0 \right\|_{L^q} \leq C [h_d(t)]^{\frac{1}{p}-\frac{1}{q}+\frac{1}{d}}e^{-\beta t} \left\| \theta_0\right\|_{L^p\cap L^2}    
\end{equation*}
for all  $u_0 \in (L^p\cap L^2)({\mathbf M};\Gamma(T{\mathbf M}))$ and vector field $\theta_0 \in (L^p\cap L^2)({\mathbf M}; \mathbb R)$. 
\end{proof}

Throughout the rest of the paper, we suppose that
$2\leq d<p<\hat p< s$ and $r>2$ such that $\dfrac{1}{2}=\dfrac{1}{p}+\dfrac{1}{\hat p}$ and $\dfrac{1}{r}=\dfrac{1}{p}+\dfrac{1}{s}$. To establish the well-posedness of mild solutions on the half time-line to system \eqref{AbstractE}, we present the following phase spaces
 \begin{eqnarray*}
 	\mathscr{X}\!\!\!\!&=&\!\!\!\!\left\lbrace u\in C_b(\r_+, (L^p\cap L^2)(\mathbf M;\Gamma(T\mathbf{M}))),  \nabla u\in  C_b(\r_+, (L^{\hat p}\cap L^s) (\mathbf M;\Gamma(T\mathbf{M})))\hbox{ such that } \right.\cr
 	&&\left. \!\!\!\! \sup\limits_{t >0}\{ \|u(t)\|_{L^p\cap L^2} + [h_d(t)]^{-(\frac{1}{p}-\frac{1}{\hat p}+\frac{1}{d})} \norm{\nabla u(t)}_{L^{\hat p}}+ [h_d(t)]^{-(\frac{1}{p}-\frac{1}{s}+\frac{1}{d})} \norm{\nabla u(t)}_{L^{s}}\} <\infty \right\rbrace 
 \end{eqnarray*}
 equipped with the norm 
 \begin{equation}\label{space1}
 	\norm{u}_{\mathscr{X}} = \sup_{t>0} \left\{ \|u(t)\|_{L^p\cap L^2} + [h_d(t)]^{-(\frac{1}{p}-\frac{1}{\hat p}+\frac{1}{d})} \norm{\nabla u(t)}_{L^{\hat p}}+ [h_d(t)]^{-(\frac{1}{p}-\frac{1}{s}+\frac{1}{d})} \norm{\nabla u(t)}_{L^{s}}\right\};
 \end{equation}
 and
  \begin{eqnarray*}
 	\mathscr{S}\!\!\!\!&=&\!\!\!\!\left\{\theta\in C_b(\r_+, (L^p\cap L^2)(\mathbf M;\r)), \, \nabla \theta\in  C_b(\r_+, (L^{\hat p} \cap L^s)(\mathbf M;\r))\hbox{ such that }\right.\cr
 	&&\left. \!\!\!\! \sup\limits_{t >0}\{ \|\theta(t)\|_{L^p\cap L^2} + [h_d(t)]^{-(\frac{1}{p}-\frac{1}{\hat p}+\frac{1}{d})} \norm{\nabla \theta(t)}_{L^{\hat p}}+ [h_d(t)]^{-(\frac{1}{p}-\frac{1}{s}+\frac{1}{d})} \norm{\nabla \theta(t)}_{L^{s}}\} <\infty \right\}
 \end{eqnarray*}
 equipped with the norm 
 \begin{equation}\label{space1}
 	\norm{\theta}_{\mathscr{S}} = \sup_{t>0} \left\{ \|\theta(t)\|_{L^p\cap L^2} + [h_d(t)]^{-(\frac{1}{p}-\frac{1}{\hat p}+\frac{1}{d})} \norm{\nabla \theta(t)}_{L^{\hat p}}+ [h_d(t)]^{-(\frac{1}{p}-\frac{1}{s}+\frac{1}{d})} \norm{\nabla \theta(t)}_{L^{s}}\right\}.
 \end{equation} 
We define and use the norm of the product space   $ \mathscr X \times \mathscr S$ by
 $\norm{\begin{bmatrix}
		u\\ \theta
\end{bmatrix}}_{\mathscr X\times \mathscr S}=\sup\limits_{t >0} \norm{\begin{bmatrix}
		u\\ \theta
	\end{bmatrix}(t)}^\blacklozenge$, where
 \begin{eqnarray*}
 	\norm{\begin{bmatrix}
 			u\\ \theta
 		\end{bmatrix}(t)}^\blacklozenge =\norm{\begin{bmatrix}
 			u\\ \theta
 		\end{bmatrix}(t)}_{L^p\cap L^2}  +& [h_d(t)]^{-(\frac{1}{p}-\frac{1}{\hat p}+\frac{1}{d})} \norm{\begin{bmatrix}
 			\nabla u\\ \nabla\theta
 		\end{bmatrix}(t)}_{L^{\hat p}} + [h_d(t)]^{-(\frac{1}{p}-\frac{1}{s}+\frac{1}{d})} \norm{\begin{bmatrix}
 			\nabla u\\ \nabla\theta
 		\end{bmatrix}(t)}_{L^{s}}. 
 \end{eqnarray*}
\begin{remark}
Note that, there are some time weights and gradient terms in the phase space $\mathscr X\times \mathscr S$. This comes from the non-commutation of $\mathbb{P}$ and $L+G=\overrightarrow{\Delta}+r+G$ (hence non-commutation of $\mathbb{P}$ and $e^{t(L+G)}$) in Remark \ref{Rem1}. This can be considered as an extension of the previous work on hyperbolic spaces \cite{XuanTrung2022}, where $\mathbb{P}$ is commutated with the vectorial heat semigroup $e^{tL}$.
\end{remark} 

Now, we can point out the existence of the bounded mild solution to the linear equation \eqref{LinearE} in the following theorem.
\begin{theorem}\label{thm1}
Let $(\mathbf{M},g)$ be a $d$-dimensional non-compact manifold (with dimension $d \geqslant 2$) satisfying Assmption \ref{Ass}. 
 Suppose that $\begin{bmatrix}
 	u_0\\ \theta_0
 \end{bmatrix}\in (L^p\cap L^2)(\mathbf M;\Gamma(T\mathbf{M}))\times (L^p\cap L^2)(\mathbf M;\r)$,  $\eta \in \mathscr S,$ the external forces $h\in \mathscr H :=C_b(\r_+, (L^{\hat p}\cap L^s)({\mathbf M};\Gamma(T{\mathbf M})))$, 
 $F\in \mathscr F:=C_b(\r_+, (L^p\cap L^2)(\mathbf M;\Gamma(T\mathbf{M})))$, $f\in \mathscr O:=C_b(\r_+, (L^p\cap L^2)(\mathbf M;\r))$. Then, the linear integral equation \eqref{LinearE} has a unique mild solution $\begin{bmatrix}
		u\\ \theta
	\end{bmatrix}\in \mathscr{X}\times \mathscr S$
  satisfying 
\begin{equation}\label{Bounded}
 \left\| \begin{bmatrix}
u\\
\theta
\end{bmatrix} \right\|_{\mathscr{X}\times \mathscr S} \leq \widetilde{C} \left\| \begin{bmatrix}
u_0\\
\theta_0
\end{bmatrix}\right\|_{L^p\cap L^2} + M \left\| h \right\|_{\mathscr H} \norm{\eta}_{\mathscr S} + N \norm{\begin{bmatrix}
F\\ f
\end{bmatrix}}_{\mathscr F\times \mathscr O},  
\end{equation}
 where the positive constants $M$ and $N$ are independent to  $h$, $\eta$, $F$ and $f$.
\end{theorem}

\begin{proof}
By Assumption \ref{Assum} and interpolation inequality (see inequality (2.7) in [30, Lemma 2.1]), we obtain that $h\in C_b(\r_+,L^q(\Gamma(T{\bf M})))$ for all $\dfrac{d}{2}<q\leq \infty$. Therefore, we have $h \in C_b(\r_+, (L^{\hat p}\cap L^s)({\bf M}; \Gamma(T{\bf M})))$ for $d<\hat p<s$. Now,  we prove the boundedness of the righ-hand side of integral equation \eqref{LinearE} which leads to the existence of linear system.

 {\bf First estimation.} By using Lemma \ref{estimates} and the boundedness of the operator $\mathbb{P}$ (see more details in example \cite{Loho}), we give the below estimate for $\left\| \begin{bmatrix}
	u(t)\\
	\theta(t)
\end{bmatrix} \right\|_{L^2}$. 
\begin{eqnarray}\label{ine-2}
\left\| \begin{bmatrix}
	u(t)\\
	\theta(t)
\end{bmatrix} \right\|_{L^2}
&\leq& \norm{e^{-t\A} \begin{bmatrix}
	u_0\\
	\theta_0
	\end{bmatrix}}_{L^2} + \norm{\mathcal H_h(\eta)(t)}_{L^2}  +\norm{\mathcal{F}\left( \begin{bmatrix}
	F\\f
	\end{bmatrix}\right) (t)}_{L^2} \cr
&\leq& C\norm{\begin{bmatrix}
		u_0\\  \theta_0
	\end{bmatrix}}_{L^2}    +\int_0^t \norm{e^{-(t-\tau)\cal{A}}\begin{bmatrix}\mathbb{P}[h\eta](\tau)\\0\end{bmatrix}}_{L^2} d\tau +\int_0^t \norm{e^{-(t-\tau)\cal{A}}  \begin{bmatrix}
	\mathbb{P}[F]\\f
	\end{bmatrix}(\tau)}_{L^2} d\tau\cr
&\leq& C\norm{\begin{bmatrix}
		u_0\\  \theta_0
	\end{bmatrix}}_{L^2} + C\int_0^t  e^{-\beta(t-\tau)}\norm{\begin{bmatrix}[h\eta](\tau)\\0\end{bmatrix}}_{L^2} d\tau  + C\int_0^t  e^{-\beta(t-\tau)} \left\|\begin{bmatrix}
	F\\f
	\end{bmatrix} (\tau) \right\|_{L^2} d\tau
\cr
&\leq& C\norm{\begin{bmatrix}
		u_0\\  \theta_0
\end{bmatrix}}_{L^2} + C\int_0^t  e^{-\beta(t-\tau)}\norm{h(\tau)}_{\hat p}\norm{\eta(\tau)}_{L^p} d\tau  + C\int_0^t  e^{-\beta(t-\tau)} \left\|\begin{bmatrix}
	F\\f
\end{bmatrix}(\tau) \right\|_{L^2} d\tau
\cr
&\leq&  C\norm{\begin{bmatrix}
			u_0\\  \theta_0
	\end{bmatrix}}_{L^2}   + C\dfrac{1}{\beta} \left( \norm{h}_{\mathscr H}\norm{\eta}_{\mathscr S}  +\norm{\begin{bmatrix}
	F\\ f
\end{bmatrix}}_{\mathscr F\times \mathscr O}\right).
	\end{eqnarray}

{\bf Second estimation.} The term $\left\| \begin{bmatrix}
	u(t)\\
	\theta(t)
\end{bmatrix} \right\|_p$ can be estimated
by using again Lemma \ref{estimates} and the boundedness of the operator $\mathbb{P}$ as follows
\begin{eqnarray}\label{ine-p}
&& \left\| \begin{bmatrix}
 	u(t)\\
 	\theta(t)
 \end{bmatrix} \right\|_{L^p}
\leq \norm{e^{-t\A} \begin{bmatrix}
 		u_0\\
 		\theta_0
 \end{bmatrix}}_{L^p} + \norm{\mathcal H_h(\eta)(t)}_{L^p}  +\norm{\mathcal{F}\left( \begin{bmatrix}
 		F\\f
 	\end{bmatrix}\right) (t)}_{L^p} \cr
 &\leq& C\norm{\begin{bmatrix}
 		u_0\\  \theta_0
 \end{bmatrix}}_{L^p\cap L^2}    +\int_0^t \norm{e^{-(t-\tau)\cal{A}}\begin{bmatrix}\mathbb{P}[h\eta](\tau)\\0\end{bmatrix}}_{L^p} d\tau +\int_0^t \norm{e^{-(t-\tau)\cal{A}}  \begin{bmatrix}
 		\mathbb{P}[F]\\f
 	\end{bmatrix}(\tau)}_{L^p} d\tau\cr
 &\leq& C\norm{\begin{bmatrix}
 		u_0\\  \theta_0
 \end{bmatrix}}_{L^p\cap L^2} + C\int_0^t\left[h_d(t-\tau)\right]^{\frac{1}{r}-\frac{1}{p}}    e^{-\beta(t-\tau)}\norm{\begin{bmatrix}[h\eta](\tau)\\0\end{bmatrix}}_{L^r\cap L^2} d\tau + C\int_0^t  e^{-\beta(t-\tau)} \left\|\begin{bmatrix}
 	F\\f
 \end{bmatrix} (\tau) \right\|_{L^p\cap L^2} d\tau
 \cr
 &\leq& C\norm{\begin{bmatrix}
 		u_0\\  \theta_0
 \end{bmatrix}}_{L^p\cap L^2} + C\int_0^t \left[h_d(t-\tau)\right]^{\frac{1}{s}}    e^{-\beta(t-\tau)} \norm{h(\tau)}_{L^s\cap L^{\hat p}}\norm{\eta(\tau)}_p  d\tau + C\int_0^t  e^{-\beta(t-\tau)} \left\|\begin{bmatrix}
 	F\\f
 \end{bmatrix}(\tau) \right\|_{L^p\cap L^2} d\tau
 \cr
 &
 \leq&C\norm{\begin{bmatrix}
 		u_0\\  \theta_0
 \end{bmatrix}}_{L^p\cap L^2}   + C\left[\beta ^{\frac{d}{2s}-1} {\bf \Gamma} (1-\frac{d}{2s})+ \dfrac{1}{\beta} \right]  \norm{h}_{\mathscr H}\norm{\eta}_{\mathscr S}  +\dfrac{C}{\beta}\norm{\begin{bmatrix}
 		F\\ f
 \end{bmatrix}}_{\mathscr F\times \mathscr O}\cr
&\leq&C\norm{\begin{bmatrix}
		u_0\\  \theta_0
\end{bmatrix}}_{L^p\cap L^2}   + M_{p_2} \norm{h}_{\mathscr H}\norm{\eta}_{\mathscr S}  +N_{p_2}\norm{\begin{bmatrix}
		F\\ f
\end{bmatrix}}_{\mathscr F\times \mathscr O}.
 \end{eqnarray}
Here $M_{p_2}=C\beta ^{\frac{d}{2s}-1} {\bf \Gamma} (1-\frac{d}{2s})+ \dfrac{2C}{\beta}$ and $N_{p_2}=\dfrac{2C}{\beta} $.

{\bf Third estimation.} Similar as above, the boundedness of the third term is given as
\begin{eqnarray}\label{ine-p^}
&&[h_d(t)]^{-(\frac{1}{p}-\frac{1}{\hat p}+\frac{1}{d})}  \left\| \begin{bmatrix}
		\nabla u(t)\\  \nabla \theta(t)
	\end{bmatrix} \right\|_{L^{\hat p}}\leq C\norm{\begin{bmatrix}
			u_0\\  \theta_0
	\end{bmatrix}}_{L^p\cap L^2}    +[h_d(t)]^{-(\frac{1}{p}-\frac{1}{\hat p}+\frac{1}{d})} \int_0^t \norm{\nabla e^{-(t-\tau)\cal{A}}\begin{bmatrix}\mathbb{P}[h\eta](\tau)\\0\end{bmatrix}}_{L^{\hat p}} d\tau  \cr
&&\hspace{6cm}+[h_d(t)]^{-(\frac{1}{p}-\frac{1}{\hat p}+\frac{1}{d})} \int_0^t \norm{\nabla e^{-(t-\tau)\cal{A}}  \begin{bmatrix}
			\mathbb{P}[F]\\f
		\end{bmatrix}(\tau)}_{L^{\hat p}} d\tau\cr
&\leq& C\norm{\begin{bmatrix}
			u_0\\  \theta_0
	\end{bmatrix}}_{L^p\cap L^2} + C\int_0^t [h_d(t)]^{-(\frac{1}{p}-\frac{1}{\hat p}+\frac{1}{d})} \left[h_d(t-\tau)\right]^{\frac{1}{r}-\frac{1}{\hat p}+\frac{1}{d}}    e^{-\beta(t-\tau)}\norm{\begin{bmatrix}[h\eta](\tau)\\0\end{bmatrix}}_{L^r\cap L^2} d\tau \cr
&& + C\int_0^t [h_d(t)]^{-(\frac{1}{p}-\frac{1}{\hat p}+\frac{1}{d})}[h_d(t-\tau)]^{\frac{1}{p}-\frac{1}{\hat p}+\frac{1}{d}} e^{-\beta(t-\tau)} \left\|\begin{bmatrix}
		F\\f
	\end{bmatrix} (\tau) \right\|_{L^p\cap L^2} d\tau
	\cr
&\leq& C\norm{\begin{bmatrix}
		u_0\\  \theta_0
\end{bmatrix}}_{L^p\cap L^2} 
+ C\int_0^t [h_d(t)]^{-(\frac{1}{p}-\frac{1}{\hat p}+\frac{1}{d})} \left[h_d(t-\tau)\right]^{\frac{1}{r}-\frac{1}{\hat p}+\frac{1}{d}}    e^{-\beta(t-\tau)}\left(\norm{h(\tau)}_{s} + \norm{h(\tau)}_{\hat p}\right)\norm{\eta(\tau)}_p d\tau 
\cr&&\!\!\! + C\int_0^t [h_d(t)]^{-(\frac{1}{p}-\frac{1}{\hat p}+\frac{1}{d})}[h_d(t-\tau)]^{\frac{1}{p}-\frac{1}{\hat p}+\frac{1}{d}} e^{-\beta(t-\tau)} \left\|\begin{bmatrix}
	F\\f
\end{bmatrix} (\tau) \right\|_{L^p\cap L^2} d\tau
\cr 
&\leq& C\norm{\begin{bmatrix}
			u_0\\  \theta_0
	\end{bmatrix}}_{L^p\cap L^2}   + M_1(t)  \norm{h}_{\mathscr H}\norm{\eta}_{\mathscr S}  +N_1(t) \norm{\begin{bmatrix}
			F\\ f
	\end{bmatrix}}_{\mathscr F\times \mathscr O}
\cr 
&\leq& C\norm{\begin{bmatrix}
		u_0\\  \theta_0
\end{bmatrix}}_{L^p\cap L^2}   + \hat M  \norm{h}_{\mathscr H}\norm{\eta}_{\mathscr S}  +\hat N \norm{\begin{bmatrix}
		F\\ f
\end{bmatrix}}_{\mathscr F\times \mathscr O},
\end{eqnarray}
where 
\begin{eqnarray*}
&& M_1(t)=	C\int_0^t [h_d(t)]^{-(\frac{1}{p}-\frac{1}{\hat p}+\frac{1}{d})} \left[h_d(t-\tau)\right]^{\frac{1}{r}-\frac{1}{\hat p}+\frac{1}{d}}    e^{-\beta(t-\tau)}d\tau \leq \hat M<+\infty,
	\cr&& N_1(t)= C\int_0^t [h_d(t)]^{-(\frac{1}{p}-\frac{1}{\hat p}+\frac{1}{d})}[h_d(t-\tau)]^{\frac{1}{p}-\frac{1}{\hat p}+\frac{1}{d}} e^{-\beta(t-\tau)}  d\tau\leq \hat N<+\infty.
\end{eqnarray*} 
 It is the fact that $N_1(t)$ is estimated in the same way of $M_1(t)$. Therefore, we merely estimate for $M_1(t)$ as follows.
\begin{eqnarray*}
	 M_1(t)=	C\int_0^t [h_d(t)]^{-(\frac{1}{p}-\frac{1}{\hat p}+\frac{1}{d})} \left[h_d(t-\tau)\right]^{\frac{1}{r}-\frac{1}{\hat p}+\frac{1}{d}}    e^{-\beta(t-\tau)}d\tau.   
\end{eqnarray*} 
 It not hard to see that $0<[h_d(t)]^{-(\frac{1}{p}-\frac{1}{\hat p}+\frac{1}{d})}<1$ for all $t>0$. Therefore, we have that
\begin{eqnarray*}
	M_1(t)\!\!\!\!\!\!&&=	C\int_0^t \left[h_d(t-\tau)\right]^{\frac{1}{r}-\frac{1}{\hat p}+\frac{1}{d}}    e^{-\beta(t-\tau)}d\tau \cr &&\leq	C\int_0^t \left(  (t-\tau)^{-\frac{d}{2}(\frac{1}{r}-\frac{1}{\hat p}+\frac{1}{d})}+1\right)  e^{-\beta(t-\tau)}d\tau \leq C \beta^{(\gamma_{1}-1)}{\bf \Gamma}(1-\gamma_{1})+\dfrac{C}{\beta}:=\hat M,  
\end{eqnarray*}  
where  $0<\gamma_{1}=\frac{d}{2}(\frac{1}{r}-\frac{1}{\hat p}+\frac{1}{d})<1$.

{\bf Fourth estimation.} The fourth term can be estimated as follows
\begin{eqnarray}\label{ine-s}
	&&[h_d(t)]^{-(\frac{1}{p}-\frac{1}{s}+\frac{1}{d})}  \left\| \begin{bmatrix}
		\nabla u(t)\\  \nabla \theta(t)
	\end{bmatrix} \right\|_{L^{s}} \leq C\norm{\begin{bmatrix}
			u_0\\  \theta_0
	\end{bmatrix}}_{L^p\cap L^2}    +[h_d(t)]^{-(\frac{1}{p}-\frac{1}{s}+\frac{1}{d})}\int_0^t \norm{\nabla e^{-(t-\tau)\cal{A}}\begin{bmatrix}\mathbb{P}[h\eta](\tau)\\0\end{bmatrix}}_{L^s} d\tau \cr 
&&\hspace{6cm} +[h_d(t)]^{-(\frac{1}{p}-\frac{1}{s}+\frac{1}{d})}\int_0^t \norm{\nabla e^{-(t-\tau)\cal{A}}  \begin{bmatrix}
			\mathbb{P}[F]\\f
		\end{bmatrix}(\tau)}_{L^s} d\tau\cr
&\leq& C\norm{\begin{bmatrix}
			u_0\\  \theta_0
	\end{bmatrix}}_{L^p\cap L^2} + C\int_0^t [h_d(t)]^{-(\frac{1}{p}-\frac{1}{s}+\frac{1}{d})} \left[h_d(t-\tau)\right]^{\frac{1}{r}-\frac{1}{s}+\frac{1}{d}}    e^{-\beta(t-\tau)}\norm{\begin{bmatrix}[h\eta](\tau)\\0\end{bmatrix}}_{L^r\cap L^2} d\tau \cr
&& + C\int_0^t [h_d(t)]^{-(\frac{1}{p}-\frac{1}{s}+\frac{1}{d})}[h_d(t-\tau)]^{\frac{1}{p}-\frac{1}{s}+\frac{1}{d}} e^{-\beta(t-\tau)} \left\|\begin{bmatrix}
		F\\f
	\end{bmatrix} (\tau) \right\|_{L^p\cap L^2} d\tau
	\cr
&\leq& C\norm{\begin{bmatrix}
			u_0\\  \theta_0
	\end{bmatrix}}_{L^p\cap L^2} 
+ C\int_0^t [h_d(t)]^{-(\frac{1}{p}-\frac{1}{s}+\frac{1}{d})} \left[h_d(t-\tau)\right]^{\frac{1}{r}-\frac{1}{s}+\frac{1}{d}}    e^{-\beta(t-\tau)}\left(\norm{h(\tau)}_{s} + \norm{h(\tau)}_{\hat p}\right)\norm{\eta(\tau)}_p d\tau
	\cr
&&+ C\int_0^t [h_d(t)]^{-(\frac{1}{p}-\frac{1}{s}+\frac{1}{d})}[h_d(t-\tau)]^{\frac{1}{p}-\frac{1}{s}+\frac{1}{d}} e^{-\beta(t-\tau)} \left\|\begin{bmatrix}
		F\\f
	\end{bmatrix} (\tau) \right\|_{L^p\cap L^2} d\tau
	\cr 
&\leq& C\norm{\begin{bmatrix}
			u_0\\  \theta_0
	\end{bmatrix}}_{L^p\cap L^2}   + M_2(t)  \norm{h}_{\mathscr H}\norm{\eta}_{\mathscr S}  +N_2(t) \norm{\begin{bmatrix}
			F\\ f
	\end{bmatrix}}_{\mathscr F\times \mathscr O}\cr 
&\leq& C\norm{\begin{bmatrix}
			u_0\\  \theta_0
	\end{bmatrix}}_{L^p\cap L^2}   + M_s \norm{h}_{\mathscr H}\norm{\eta}_{\mathscr S}  + N_s \norm{\begin{bmatrix}
			F\\ f
	\end{bmatrix}}_{\mathscr F\times \mathscr O},
\end{eqnarray}
where 
\begin{eqnarray*}
	&& M_2(t)=	C\int_0^t [h_d(t)]^{-(\frac{1}{p}-\frac{1}{s}+\frac{1}{d})} \left[h_d(t-\tau)\right]^{\frac{1}{r}-\frac{1}{s}+\frac{1}{d}}   e^{-\beta(t-\tau)}d\tau \leq M_s<+\infty,
	\cr&& N_2(t)= C\int_0^t [h_d(t)]^{-(\frac{1}{p}-\frac{1}{s}+\frac{1}{d})}[h_d(t-\tau)]^{\frac{1}{p}-\frac{1}{s}+\frac{1}{d}} e^{-\beta(t-\tau)}  d\tau\leq N_s<+\infty.
\end{eqnarray*} 
It is the fact that $M_2(t), N_2(t)$ are estimated in the same way of $M_1(t)$. Therefore, we merely estimate for $M_2(t)$ as following.
\begin{eqnarray*}
M_2(t)\!\!\!\!\!\!\!&&=	C\int_0^t [h_d(t)]^{-(\frac{1}{p}-\frac{1}{s}+\frac{1}{d})} \left[h_d(t-\tau)\right]^{\frac{1}{r}-\frac{1}{s}+\frac{1}{d}}   e^{-\beta(t-\tau)}d\tau   
\cr&&\leq	C\int_0^t \left(  (t-\tau)^{-\frac{d}{2}(\frac{1}{r}-\frac{1}{s}+\frac{1}{d})}+1\right)  e^{-\beta(t-\tau)}d\tau \leq C\beta^{(\gamma_{2}-1)}{\bf \Gamma}(1-\gamma_{2})+\dfrac{C}{\beta}.  
\end{eqnarray*}  
Here  $0<\gamma_{2}=\frac{d}{2}(\frac{1}{r}-\frac{1}{s} +\frac{1}{d})<1$ and $M_s= C\beta^{(\gamma_{2}-1)}{\bf \Gamma}(1-\gamma_{2})+\dfrac{C}{\beta}$. 
 
Finally, setting $ M=M_{p_2}+\hat{M}+M_s$ and $N=N_{p_2}+\hat{N}+N_s$ and combining the inequalities \eqref{ine-2}, \eqref{ine-p},  \eqref{ine-p^} and \eqref{ine-s}, we obtain the boundedness \eqref{Bounded}. Our proof is completed.
\end{proof}

\subsection{Bilinear estimates and well-posedness for Boussinesq systems}
To prove the well-posedness of mild solutions for semilinear equation \eqref{MatrixEq}, we first obtain the bilinear estimates for operator $\mathscr{B}(\cdot,\cdot)$ defining by \eqref{Bilinear}.
\begin{lemma}\label{BiE}
Let $({\bf M},g)$ be a $d$-dimensional non-compact manifold (with dimension $d\geqslant 2$) satisfying Assumption \ref{Ass}. There exists a universal constant $K>0$ such that

(i) for all $t>0$, 
\begin{equation}\label{Bilinear-2}
		\norm{\mathscr B\left( \begin{bmatrix}
				u\\ \zeta   
			\end{bmatrix}, \begin{bmatrix}
				v\\ \vartheta   
			\end{bmatrix}  \right)(t)}_{2} \leq K \norm{\begin{bmatrix}
			u\\
			\zeta
		\end{bmatrix}}_{\mathscr X\times\mathscr S}\norm{\begin{bmatrix}
		v\\ \vartheta
	\end{bmatrix}}_{\mathscr X\times\mathscr S};  
	\end{equation}
	
	(ii) for all $t>0$, \begin{equation}\label{Bilinear-p}
		\norm{\mathscr B\left( \begin{bmatrix}
				u\\ \zeta   
			\end{bmatrix}, \begin{bmatrix}
				v\\ \vartheta   
			\end{bmatrix}  \right)(t)}_{p} \leq K \norm{\begin{bmatrix}
				u\\
				\zeta
		\end{bmatrix}}_{\mathscr X\times\mathscr S}\norm{\begin{bmatrix}
				v\\ \vartheta
	\end{bmatrix}}_{\mathscr X\times\mathscr S};  
		\end{equation}
	
(iii) for all $t>0$,
\begin{equation}\label{Bilinear-p^}
	[h_d(t)]^{-(\frac{1}{p}-\frac{1}{\hat p}+\frac{1}{d})}\norm{\nabla \mathscr B\left( \begin{bmatrix}
			u\\ \zeta   
		\end{bmatrix}, \begin{bmatrix}
			v\\ \vartheta   
		\end{bmatrix}  \right)(t)}_{\hat p} \leq K \norm{\begin{bmatrix}
			u\\
			\zeta
	\end{bmatrix}}_{\mathscr X\times\mathscr S}\norm{\begin{bmatrix}
			v\\ \vartheta
	\end{bmatrix}}_{\mathscr X\times\mathscr S};  
	\end{equation}
	
(iv) and for all $t>0$,
\begin{equation}\label{Bilinear-s}
	[h_d(t)]^{-(\frac{1}{p}-\frac{1}{s}+\frac{1}{d})}\norm{\nabla \mathscr B\left( \begin{bmatrix}
			u\\ \zeta   
		\end{bmatrix}, \begin{bmatrix}
			v\\ \vartheta   
		\end{bmatrix}  \right)(t)}_{s} \leq K \norm{\begin{bmatrix}
			u\\
			\zeta
	\end{bmatrix}}_{\mathscr X\times\mathscr S}\norm{\begin{bmatrix}
			v\\ \vartheta
	\end{bmatrix}}_{\mathscr X\times\mathscr S},
\end{equation}
where the positive constant $K$ does not dependent on funtions $u,v,\zeta,\vartheta$.
\end{lemma}
\begin{proof}
Using the boundedness of Kodaira-Hodge operator $\mathbb{P}$ and the $L^p-L^q$ smoothing estimates in Lemma \ref{estimates}$(ii)$ and H\"older's inequality, we obtain that
\begin{eqnarray}\label{Es-Bilinear-2}
&&\norm{\mathscr B\left( \begin{bmatrix}
		u\\ \zeta   
	\end{bmatrix}, \begin{bmatrix}
		v\\ \vartheta   
	\end{bmatrix}  \right)(t)}_2
\leq \int_0^t  \norm{e^{-(t-\tau)\A} \begin{bmatrix}
		\p [(u \cdot\nabla)v](\tau)\\
		[(u \cdot\nabla)\vartheta](\tau)
\end{bmatrix}}_2 d\tau  \cr
	&\leq&  2C\int_0^t  e^{-\beta(t-\tau)}\norm{\begin{bmatrix}
		 [(u \cdot\nabla)v](\tau)\\
		[(u \cdot\nabla)\vartheta](\tau)
	\end{bmatrix} }_2d\tau  \leq  2C\int_0^t  e^{-\beta(t-\tau)}\norm{u(\tau)}_{ p}\norm{\begin{bmatrix}
		\nabla v(\tau)\\
		\nabla\vartheta(\tau)
\end{bmatrix} }_{\hat p} d\tau  \cr
	&\leq&  2C\norm{\begin{bmatrix}
			u\\   \zeta
	\end{bmatrix}}_{\mathscr X\times\mathscr S}\norm{\begin{bmatrix}
			v\\  \vartheta
	\end{bmatrix}}_{\mathscr X\times\mathscr S}  \int_0^t  [h_d(\tau)]^{\frac{1}{p}-\frac{1}{\hat p}+\frac{1}{d} } e^{-\beta(t-\tau)}d\tau\cr
	&\leq& K_2\norm{\begin{bmatrix}
			u\\  \zeta
	\end{bmatrix}}_{\mathscr X\times\mathscr S}\norm{\begin{bmatrix}
			v\\
			\vartheta
	\end{bmatrix}}_{\mathscr X\times\mathscr S},
\end{eqnarray}
where  
$$G_2(t)=2C\int_0^t  [h_d(\tau)]^{\frac{1}{p}-\frac{1}{\hat p}+\frac{1}{d}} e^{-\beta(t-\tau)}d\tau\leq K_2<+\infty.$$
Indeed, for $0<t\leq 1$ one has $G_2(t)\leq 2C\int_0^1  (\tau)^{-\frac{d}{2}(\frac{1}{p}-\frac{1}{\hat p}+\frac{1}{d})} d\tau=\dfrac{2C}{\frac{1}{2}-\frac{d}{2p}+\frac{d}{2\hat p}}.$ For the case $t>1$, we have
\begin{eqnarray*}
G_2(t)\leq 2C\int_0^1  (\tau)^{-\frac{d}{2}(\frac{1}{p}-\frac{1}{\hat p}+\frac{1}{d})} d\tau+2C\int_1^t   e^{-\beta(t-\tau)}d\tau\leq \dfrac{2C}{\frac{1}{2}-\frac{d}{2p}+\frac{d}{2\hat p}}+\dfrac{2C}{\beta}:=K_2.
\end{eqnarray*}

Similar as above, we can estimate $\norm{\mathscr B\left( \begin{bmatrix}
		u\\ \zeta   
	\end{bmatrix}, \begin{bmatrix}
		v\\ \vartheta  
	\end{bmatrix}  \right)(t)}_p$ by
\begin{eqnarray}\label{Es-Bilinear-p}
	&&\norm{\mathscr B\left( \begin{bmatrix}
			u\\ \zeta   
		\end{bmatrix}, \begin{bmatrix}
			v\\ \vartheta   
		\end{bmatrix}  \right)(t)}_p
\leq \int_0^t  \norm{e^{-(t-\tau)\A} \begin{bmatrix}
			\p [(u \cdot\nabla)v](\tau)\\
			[(u \cdot\nabla)\vartheta](\tau)
	\end{bmatrix}}_p d\tau  \cr
&\leq&  C\int_0^t [h_d(t-\tau)]^{\frac{1}{r}-\frac{1}{p}} e^{-\beta(t-\tau)}\norm{\begin{bmatrix}
			[(u \cdot\nabla)v](\tau)\\
			[(u \cdot\nabla)\vartheta](\tau)
	\end{bmatrix} }_{L^r\cap L^2}d\tau \cr
&\leq& C\int_0^t [h_d(t-\tau)]^{\frac{1}{s}} e^{-\beta(t-\tau)}\left[\norm{u(\tau)}_{ p}\norm{\begin{bmatrix}
			\nabla v(\tau)\\
			\nabla\vartheta(\tau)
	\end{bmatrix} }_{s}+ \norm{u(\tau)}_{ p}\norm{\begin{bmatrix}
			\nabla v(\tau)\\
			\nabla\vartheta(\tau)
	\end{bmatrix} }_{\hat p}\right] d\tau  \cr
&\leq& C\norm{\begin{bmatrix}
			u\\   \zeta
	\end{bmatrix}}_{\mathscr X\times\mathscr S}\norm{\begin{bmatrix}
			v\\  \vartheta
	\end{bmatrix}}_{\mathscr X\times\mathscr S}  \int_0^t  [h_d(t-\tau)]^{\frac{1}{s}}\left[ [h_d(\tau)]^{\frac{1}{p}-\frac{1}{s}+\frac{1}{d}} +[h_d(\tau)]^{\frac{1}{p}-\frac{1}{\hat p}+\frac{1}{d}} \right]e^{-\beta(t-\tau)} d\tau\cr
&\leq& K_p\norm{\begin{bmatrix}
			u\\  \zeta
	\end{bmatrix}}_{\mathscr X\times\mathscr S}\norm{\begin{bmatrix}
			v\\
			\vartheta
	\end{bmatrix}}_{\mathscr X\times\mathscr S},
\end{eqnarray}
where (see Subsection \ref{convergence}): 
$$G_p(t):=C\int_0^t  [h_d(t-\tau)]^{\frac{1}{s}}\left[ [h_d(\tau)]^{\frac{1}{p}-\frac{1}{s}+\frac{1}{d}} +[h_d(\tau)]^{\frac{1}{p}-\frac{1}{\hat p}+\frac{1}{d}} \right]e^{-\beta(t-\tau)}  d\tau\leq K_p<+\infty.$$

On the other hand, we have the following estimation 
\begin{eqnarray}\label{Es-Bilinear-p^}
&&[h_d(t)]^{-(\frac{1}{p}-\frac{1}{\hat p}+\frac{1}{d})}\norm{\nabla \mathscr B\left( \begin{bmatrix}
 			u\\ \zeta   
 		\end{bmatrix}, \begin{bmatrix}
 			v\\ \vartheta   
 		\end{bmatrix}  \right)(t)}_{\hat p}
\leq \int_0^t [h_d(t)]^{-(\frac{1}{p}-\frac{1}{\hat p}+\frac{1}{d})} \norm{\nabla e^{-(t-\tau)\A} \begin{bmatrix}
 			\p [(u \cdot\nabla)v](\tau)\\
 			[(u \cdot\nabla)\vartheta](\tau)
 	\end{bmatrix}}_{\hat p} d\tau  \cr
&\leq&  C\int_0^t [h_d(t)]^{-(\frac{1}{p}-\frac{1}{\hat p}+\frac{1}{d})} [h_d(t-\tau)]^{\frac{1}{r}-\frac{1}{\hat p}+\frac{1}{d}} e^{-\beta(t-\tau)}\norm{\begin{bmatrix}
 			[(u \cdot\nabla)v](\tau)\\
 			[(u \cdot\nabla)\vartheta](\tau)
 	\end{bmatrix} }_{L^r\cap L^2}d\tau \cr
&\leq& C\int_0^t [h_d(t)]^{-(\frac{1}{p}-\frac{1}{\hat p}+\frac{1}{d})} [h_d(t-\tau)]^{\frac{1}{r}-\frac{1}{\hat p}+\frac{1}{d}} e^{-\beta(t-\tau)}\cr
 	&& \hspace{5cm}\times \left[\norm{u(\tau)}_{ p}\norm{\begin{bmatrix}
 			\nabla v(\tau)\\
 			\nabla\vartheta(\tau)
 	\end{bmatrix} }_{s}+ \norm{u(\tau)}_{ p}\norm{\begin{bmatrix}
 			\nabla v(\tau)\\
 			\nabla\vartheta(\tau)
 	\end{bmatrix} }_{\hat p}\right] d\tau  \cr
&\leq& C\norm{\begin{bmatrix}
 			u\\   \zeta
 	\end{bmatrix}}_{\mathscr X\times\mathscr S}\norm{\begin{bmatrix}
 			v\\  \vartheta
 	\end{bmatrix}}_{\mathscr X\times\mathscr S}  \int_0^t [h_d(t)]^{-(\frac{1}{p}-\frac{1}{\hat p}+\frac{1}{d})} [h_d(t-\tau)]^{\frac{1}{r}-\frac{1}{\hat p}+\frac{1}{d}}  \left[ [h_d(\tau)]^{\frac{1}{p}-\frac{1}{s}+\frac{1}{d}} +[h_d(\tau)]^{\frac{1}{p}-\frac{1}{\hat p}+\frac{1}{d}} \right]e^{-\beta(t-\tau)} d\tau\cr
&\leq& K_{\hat p} \norm{\begin{bmatrix}
 			u\\  \zeta
 	\end{bmatrix}}_{\mathscr X\times\mathscr S}\norm{\begin{bmatrix}
 			v\\
 			\vartheta
 	\end{bmatrix}}_{\mathscr X\times\mathscr S}.
 \end{eqnarray}
Here, we used the fact that $ [h_d(\tau)]^{\frac{1}{p}-\frac{1}{s}+\frac{1}{d}} >[h_d(\tau)]^{\frac{1}{p}-\frac{1}{\hat p}+\frac{1}{d}} $ for all $\tau>0$, and by the samw way as the boundedness of $K_s$ we have the following boundedness
 $$G_{\hat p}(t):=2C \int_0^t [h_d(t)]^{-(\frac{1}{p}-\frac{1}{\hat p}+\frac{1}{d})} [h_d(t-\tau)]^{\frac{1}{r}-\frac{1}{\hat p}+\frac{1}{d}}  [h_d(\tau)]^{\frac{1}{p}-\frac{1}{s}+\frac{1}{d}}  e^{-\beta(t-\tau)} d\tau\leq K_{\hat p}<+\infty.$$

Finally, the boundedness of $[h_d(t)]^{-(\frac{1}{p}-\frac{1}{s}+\frac{1}{d})}\norm{\nabla \mathscr B\left( \begin{bmatrix}
		u\\ \zeta   
	\end{bmatrix}, \begin{bmatrix}
		v\\ \vartheta   
	\end{bmatrix}  \right)(t)}_{s}$ is given as follows
\begin{eqnarray}\label{Es-Bilinear-s}
	&&[h_d(t)]^{-(\frac{1}{p}-\frac{1}{s}+\frac{1}{d})}\norm{\nabla \mathscr B\left( \begin{bmatrix}
			u\\ \zeta   
		\end{bmatrix}, \begin{bmatrix}
			v\\ \vartheta   
		\end{bmatrix}  \right)(t)}_{s}
\leq \int_0^t [h_d(t)]^{-(\frac{1}{p}-\frac{1}{s}+\frac{1}{d})} \norm{\nabla e^{-(t-\tau)\A} \begin{bmatrix}
			\p [(u \cdot\nabla)v](\tau)\\
			[(u \cdot\nabla)\vartheta](\tau)
	\end{bmatrix}}_s d\tau  \cr
&\leq&  C\int_0^t [h_d(t)]^{-(\frac{1}{p}-\frac{1}{s}+\frac{1}{d})} [h_d(t-\tau)]^{\frac{1}{r}-\frac{1}{s}+\frac{1}{d}} e^{-\beta(t-\tau)}\norm{\begin{bmatrix}
			[(u \cdot\nabla)v](\tau)\\
			[(u \cdot\nabla)\vartheta](\tau)
	\end{bmatrix} }_{L^r\cap L^2}d\tau \cr
	&\leq& C\int_0^t [h_d(t)]^{-(\frac{1}{p}-\frac{1}{s}+\frac{1}{d})} [h_d(t-\tau)]^{\frac{1}{r}-\frac{1}{s}+\frac{1}{d}} e^{-\beta(t-\tau)} \left[\norm{u(\tau)}_{ p}\norm{\begin{bmatrix}
			\nabla v(\tau)\\
			\nabla\vartheta(\tau)
	\end{bmatrix} }_{s}+ \norm{u(\tau)}_{ p}\norm{\begin{bmatrix}
			\nabla v(\tau)\\
			\nabla\vartheta(\tau)
	\end{bmatrix} }_{\hat p}\right] d\tau \cr
&\leq& C\norm{\begin{bmatrix}
			u\\   \zeta
	\end{bmatrix}}_{\mathscr X\times\mathscr S}\norm{\begin{bmatrix}
			v\\  \vartheta
	\end{bmatrix}}_{\mathscr X\times\mathscr S}  \int_0^t [h_d(t)]^{-(\frac{1}{p}-\frac{1}{s}+\frac{1}{d})} [h_d(t-\tau)]^{\frac{1}{r}-\frac{1}{s}+\frac{1}{d}}  \left[ [h_d(\tau)]^{\frac{1}{p}-\frac{1}{s}+\frac{1}{d}} +[h_d(\tau)]^{\frac{1}{p}-\frac{1}{\hat p}+\frac{1}{d}} \right]e^{-\beta(t-\tau)} d\tau\cr
&\leq& K_{s}\norm{\begin{bmatrix}
			u\\  \zeta
	\end{bmatrix}}_{\mathscr X\times\mathscr S}\norm{\begin{bmatrix}
			v\\
			\vartheta
	\end{bmatrix}}_{\mathscr X\times\mathscr S}.
\end{eqnarray}
Here, we used the following integral's convergence  (see Subsection \ref{convergence}):
$$G_{s}(t):=2C \int_0^t [h_d(t)]^{-(\frac{1}{p}-\frac{1}{s}+\frac{1}{d})} [h_d(t-\tau)]^{\frac{1}{r}-\frac{1}{s}+\frac{1}{d}} [h_d(\tau)]^{\frac{1}{p}-\frac{1}{s}+\frac{1}{d}}  e^{-\beta(t-\tau)} d\tau\leq K_s<+\infty.$$ 
Therefore, by setting $K=\max \{K_2, K_p, K_{\hat p}, K_s\}$, we obtain the inequalities \eqref{Bilinear-2}, \eqref{Bilinear-p}, \eqref{Bilinear-p^} and \eqref{Bilinear-s} as desired.
\end{proof}

We state and prove the well-posedness of mild solutions for semilinear equation \eqref{MatrixEq} in the following theorem.
\begin{theorem}\label{PeriodicThm}(Global-in -time mild solution). Let $({\bf M},g)$ be a $d$-dimensional noncompact manifold (with dimension $d\geq 2$) satisfying Assumption \ref{Ass} . For $p> d$, suppose that  the external forces $h\in \mathscr H:= C_b(\r_+, (L^{\hat p}\cap L^{s})({\mathbf M};\Gamma(T{\mathbf M})))$, 
$F\in \mathscr F:=C_b(\r_+, (L^{p}\cap L^{2})(\mathbf M;\Gamma(T\mathbf{M})))$, $f\in \mathscr O:=C_b(\r_+, (L^{p}\cap L^{2})(\mathbf M;\r))$. If the norms $\norm{\begin{bmatrix}
		u_0\\ \theta_0
\end{bmatrix}}_{L^p\cap L^2}$, $\norm{h}_{\mathscr H}$, $\norm{\begin{bmatrix}
		F\\ f
\end{bmatrix}}_{\mathscr F\times \mathscr O}$  are sufficiently small, then equation \eqref{MatrixEq} has one and only one bounded mild solution $(\hat{u},\hat{\theta})$ on a small ball of  
$\mathscr X\times \mathscr S$.
\end{theorem} 
\begin{proof}
In order to start, we denote 
\begin{eqnarray*}\label{bro}
\B_\rho:=\left\lbrace \begin{bmatrix}
	v\\
	\vartheta
\end{bmatrix}\in \mathscr X\times \mathscr S \hbox{ such that }   \norm{\begin{bmatrix}
	v\\
	\vartheta
	\end{bmatrix}}_{\mathscr X\times \mathscr S} \leq \rho \right\rbrace .
\end{eqnarray*}
For each $\begin{bmatrix}
	v\\
	\vartheta
\end{bmatrix}\in \B_\rho$, we consider the linear equation 
\begin{equation}\label{ns1}
\begin{bmatrix}
u(t)\\
\theta(t)
\end{bmatrix}
=e^{-t\A}%
\begin{bmatrix}
u_{0}\\
\theta_{0}%
\end{bmatrix}
+\mathscr B\left(
\begin{bmatrix}
v\\
\vartheta
\end{bmatrix}
,%
\begin{bmatrix}
v\\
\vartheta
\end{bmatrix}
\right)  (t)+\cal H_h(\vartheta)(t) + \mathcal F \left(
\begin{bmatrix}
F\\
f
\end{bmatrix}
\right)(t).
\end{equation}
Now, we apply the bilinear estimates \eqref{Bilinear-2}, \eqref{Bilinear-p}, \eqref{Bilinear-p^} and \eqref{Bilinear-s} for $\mathscr B(\cdot,\cdot)$ in Lemma \ref{BiE} and the linear estimate for $\mathcal{F}(\cdot)$ as in the proof of Theorem \ref{thm1} to obtain that for $(v,\vartheta)\in \B_\rho$  there exists a unique bounded mild solution $(u,\theta)$ to \eqref{ns1} satisfying
\begin{eqnarray}\label{ephi}
\norm{\begin{bmatrix}
u\\ \theta    
\end{bmatrix}}_{\mathscr X\times \mathscr S}&\le & \widetilde{C} \norm{ \begin{bmatrix}
u_0\\ \theta_0    
\end{bmatrix}}_{L^p\cap L^2} + 4K\norm{\begin{bmatrix}
 v \\ \vartheta   
\end{bmatrix}}^2_{\mathscr X\times \mathscr S}+ M \norm{h}_{\mathscr H}\norm{\vartheta}_{\mathscr S}+  N\norm{\begin{bmatrix}
F\\ f
\end{bmatrix}}_{\mathscr F\times \mathscr O}\cr
&\le & \widetilde{C} \norm{ \begin{bmatrix}
		u_0\\ \theta_0    
\end{bmatrix} }_{L^p\cap L^2} + 4K\rho^2+ M \norm{h}_{\mathscr H}\rho+  N\norm{\begin{bmatrix}
F\\ f
\end{bmatrix}}_{\mathscr F\times \mathscr O}.
\end{eqnarray}
Thence, we can define a map ${\bf \Phi}: \mathscr X\times \mathscr S \longrightarrow  \mathscr X\times \mathscr S$ as follows
\begin{equation}\label{defphi}
\begin{split}
{\bf\Phi} \begin{bmatrix}
  v\\\vartheta  
\end{bmatrix}&= \begin{bmatrix}
    u\\ \theta
\end{bmatrix}
\end{split}.
\end{equation}
If $\rho$, $\norm{h}_{\mathscr H}$, $\norm{\begin{bmatrix}
		u_0\\ \theta_0
\end{bmatrix} }_{L^p\cap L^2}$  and $\norm{\begin{bmatrix}
F\\ f
\end{bmatrix}}_{\mathscr F\times \mathscr O}$ are small enough, then $\norm{\begin{bmatrix}
		{\bf\Phi} \begin{bmatrix}
			v\\\vartheta  
		\end{bmatrix}
\end{bmatrix}}_{\mathscr X\times \mathscr S}<  \rho$. This means that the mapping $\Phi$ maps  $\B_\rho$ into itself.

Moreover, we have clearly that
\begin{equation}\label{defphi1}
{\bf \Phi} \left( \begin{bmatrix}
    v\\\vartheta
\end{bmatrix}\right)(t) = \begin{bmatrix}
    u_0\\ \theta_0
\end{bmatrix} +\mathscr B\left( \begin{bmatrix}
    v\\\vartheta
\end{bmatrix},\begin{bmatrix}
    v\\\vartheta
\end{bmatrix} \right)(t) + \cal H_h(\vartheta)(t) + \mathcal{F} \left( \begin{bmatrix}
    F\\ f
\end{bmatrix} \right)(t).
\end{equation}
Hence, for $(v_1,\vartheta_1),\, (v_2,\vartheta_2)\in \B_\rho$, applying again  the linear estimate for $\mathcal H(\cdot)$ as in the proof of Theorem \ref{thm1} and performing some computations, we have 
\begin{eqnarray}\label{ephi33}
&&\norm{{\bf \Phi}\left( \begin{bmatrix}
v_1\\\vartheta_1
\end{bmatrix}\right)-{\bf \Phi} \left( \begin{bmatrix}
    v_2\\\vartheta_2
\end{bmatrix}\right)}_{\mathscr X\times\mathscr S}= \norm{\mathscr B\left( \begin{bmatrix}
    v_1\\\vartheta_1
\end{bmatrix},\begin{bmatrix}
    v_1\\\vartheta_1
\end{bmatrix} \right) -\mathscr B\left( \begin{bmatrix}
    v_2\\\vartheta_2
\end{bmatrix},\begin{bmatrix}
    v_2\\\vartheta_2
\end{bmatrix} \right) + \cal H_h(\vartheta_1-\vartheta_2)}_{\mathscr X\times\mathscr S}\cr
&\leq& \sup\limits_{t >0} \norm{ \int_0^t e^{-(t-\tau)\A} \begin{bmatrix}
    \mathbb P([(v_1-v_2) \cdot\nabla ]v_1 + [v_2\cdot\nabla ](v_1-v_2))\\ 
   [(v_1-v_2) \cdot \nabla ]\vartheta_2+ [v_1\cdot \nabla](\vartheta_1-\vartheta_2) 
\end{bmatrix}(\tau)d\tau}^\blacklozenge  \cr&& +\norm{  \cal H_h(\vartheta_1-\vartheta_2)}_{\mathscr X\times\mathscr S}\cr
&\leq&  \sup\limits_{t >0} \left\lbrace 
\norm{\int_0^t e^{-(t-\tau)\A} \begin{bmatrix}
		\mathbb P([(v_1-v_2) \cdot\nabla ]v_1 + [v_2\cdot\nabla ](v_1-v_2))\\ 
		[(v_1-v_2) \cdot \nabla ]\vartheta_2+ [v_1\cdot \nabla](\vartheta_1-\vartheta_2) 
	\end{bmatrix}(\tau)d\tau}_{L^p\cap L^2}\right.\cr
&&\left. + [h_d(t)]^{-(\frac{1}{p}-\frac{1}{\hat p}+\frac{1}{d})} \norm{\int_0^t \nabla e^{-(t-\tau)\A} \begin{bmatrix}
		\mathbb P([(v_1-v_2) \cdot\nabla ]v_1 + [v_2\cdot\nabla ](v_1-v_2))\\ 
		[(v_1-v_2) \cdot \nabla ]\vartheta_2+ [v_1\cdot \nabla](\vartheta_1-\vartheta_2) 
	\end{bmatrix}(\tau)d\tau}_{L^{\hat p}} 
\right.\cr
&&\left.+ [h_d(t)]^{-(\frac{1}{p}-\frac{1}{s}+\frac{1}{d})} \norm{\int_0^t\nabla e^{-(t-\tau)\A} \begin{bmatrix}
		\mathbb P([(v_1-v_2) \cdot\nabla ]v_1 + [v_2\cdot\nabla ](v_1-v_2))\\ 
		[(v_1-v_2) \cdot \nabla ]\vartheta_2+ [v_1\cdot \nabla](\vartheta_1-\vartheta_2) 
	\end{bmatrix}(\tau)d\tau}_{L^{s}} \right\rbrace \cr
&& +\norm{  \cal H_h(\vartheta_1-\vartheta_2)}_{\mathscr X\times\mathscr S}
\cr
&\leq& \sup\limits_{t >0} \left\lbrace  8C\rho  \int_0^t [h_d(t-\tau)]^{\frac{1}{s}} [h_d(\tau)]^{\frac{1}{p}-\frac{1}{s}+\frac{1}{d}} e^{-\beta(t-\tau)}  \norm{ \begin{bmatrix}
		v_1-v_2\\\vartheta_1-\vartheta_2 \end{bmatrix}(\tau)}^\blacklozenge d\tau
\right.\cr
&&\left.+4C\rho  \int_0^t [h_d(t)]^{-(\frac{1}{p}-\frac{1}{\hat p}+\frac{1}{d})} [h_d(t-\tau)]^{\frac{1}{r}-\frac{1}{\hat  p}+\frac{1}{d}} [h_d(\tau)]^{\frac{1}{p}-\frac{1}{s}+\frac{1}{d}} e^{-\beta(t-\tau)}  \norm{ \begin{bmatrix}
	v_1-v_2\\\vartheta_1-\vartheta_2 \end{bmatrix}(\tau)}^\blacklozenge d\tau 
\right.\cr
&&\left.+4C\rho  \int_0^t [h_d(t)]^{-(\frac{1}{p}-\frac{1}{s}+\frac{1}{d})} [h_d(t-\tau)]^{\frac{1}{r}-\frac{1}{s}+\frac{1}{d}} [h_d(\tau)]^{\frac{1}{p}-\frac{1}{s}+\frac{1}{d}} e^{-\beta (t-\tau)}  \norm{ \begin{bmatrix}
	v_1-v_2\\\vartheta_1-\vartheta_2 \end{bmatrix}(\tau)}^\blacklozenge d\tau\right\rbrace 
 \cr&& +\norm{  \cal H_h(\vartheta_1-\vartheta_2)}_{\mathscr X\times\mathscr S}
 \cr
 &\leq& \sup\limits_{t >0} \left\lbrace  8C\rho  \int_0^t [h_d(t-\tau)]^{\frac{1}{s}} [h_d(\tau)]^{\frac{1}{p}-\frac{1}{s}+\frac{1}{d}} e^{-\beta(t-\tau)}  d\tau \norm{ \begin{bmatrix}
 		v_1-v_2\\\vartheta_1-\vartheta_2 \end{bmatrix}}_{\mathscr X\times\mathscr S}
\right.\cr
&&\left.+ 4C\rho  \int_0^t [h_d(t)]^{-(\frac{1}{p}-\frac{1}{\hat p}+\frac{1}{d})} [h_d(t-\tau)]^{\frac{1}{r}-\frac{1}{\hat  p}+\frac{1}{d}} [h_d(\tau)]^{\frac{1}{p}-\frac{1}{s}+\frac{1}{d}} e^{-\beta(t-\tau)} d\tau \norm{ \begin{bmatrix}
 	v_1-v_2\\\vartheta_1-\vartheta_2 \end{bmatrix}}_{\mathscr X\times\mathscr S}
\right.\cr
&&\left.+ 4C\rho \int_0^t [h_d(t)]^{-(\frac{1}{p}-\frac{1}{s}+\frac{1}{d})} [h_d(t-\tau)]^{\frac{1}{r}-\frac{1}{s}+\frac{1}{d}} [h_d(\tau)]^{\frac{1}{p}-\frac{1}{s}+\frac{1}{d}} e^{-\beta (t-\tau)}  d\tau\norm{ \begin{bmatrix}
 	v_1-v_2\\\vartheta_1-\vartheta_2 \end{bmatrix}}_{\mathscr X\times\mathscr S}\right\rbrace 
 \cr&& +\norm{  \cal H_h(\vartheta_1-\vartheta_2)}_{\mathscr X\times\mathscr S}
\cr
&\leq& \left(8CL_1\rho+4CL_2\rho+4CL_3\rho + M \norm{h}_{\mathscr H}\right) \norm{ \begin{bmatrix}
    v_1-v_2\\\vartheta_1-\vartheta_2 \end{bmatrix}}_{\mathscr X\times\mathscr S}, 
\end{eqnarray}
where, we used  
\begin{eqnarray*}
&& \int_0^t [h_d(t-\tau)]^{\frac{1}{s}} [h_d(\tau)]^{\frac{1}{p}-\frac{1}{s}+\frac{1}{d}} e^{-\beta(t-\tau)}  d\tau \leq L_1<+\infty, 
\cr&& \int_0^t [h_d(t)]^{-(\frac{1}{p}-\frac{1}{\hat p}+\frac{1}{d})} [h_d(t-\tau)]^{\frac{1}{r}-\frac{1}{\hat  p}+\frac{1}{d}} [h_d(\tau)]^{\frac{1}{p}-\frac{1}{s}+\frac{1}{d}} e^{-\beta(t-\tau)} d\tau  \leq L_2<+\infty,
\cr&&\int_0^t [h_d(t)]^{-(\frac{1}{p}-\frac{1}{s}+\frac{1}{d})} [h_d(t-\tau)]^{\frac{1}{r}-\frac{1}{s}+\frac{1}{d}} [h_d(\tau)]^{\frac{1}{p}-\frac{1}{s}+\frac{1}{d}} e^{-\beta (t-\tau)}  d\tau \leq L_3<+\infty
\end{eqnarray*} 
with the finding of bounded constans $L_1,L_2 $ and $L_3$ is similar to $K_p, K_{\hat p} $ and $ K_s$ and we omit the details.

Therefore, if $\rho$ and $\norm{h}_{\mathscr H}$ are small enough such that $8CL_1\rho+4CL_2\rho+4CL_3\rho + M \norm{h}_{\mathscr H}<1$ then the map ${\bf \Phi}$ becomes a contraction on $\B_\rho$.
By fixed point arguments, there exists a unique fixed point $(\hat{u},\hat\theta)$ of ${\bf \Phi}$, and by the definition of ${\bf \Phi}$, this fixed point $(\hat{u},\hat\theta)$ is a bounded mild solution to equation \eqref{MatrixEq}.
The uniqueness of $(\hat{u},\hat\theta)$ in the small ball $\B_\rho$ is clearly by using inequality \eqref{ephi33}.
\end{proof}

\subsection{Exponential stability}\label{S4}
In this section, we will show that how we use Gronwall's inequality to prove the exponential stability of mild solutions to integral equation \eqref{MatrixEq} obtained in Theorem \ref{PeriodicThm}. Our main theorem is as follows:
\begin{theorem}\label{stability}(Exponential stability).
	Let $({\bf M},g)$ be a $d$-dimensional noncompact manifold with $d\geqslant 2$. For $p> d$, assume that the external force $h\in C_b(\r_+, (L^{\hat p}\cap L^s)({\mathbf M};\Gamma(T{\mathbf M})))$ with the norm $\norm{h}_{\mathscr H}$ small enough. Then, the mild solution $(u,\theta)$ of the equation \eqref{MatrixEq} is exponentially stable in the sense that for any other mild solution $\begin{bmatrix}
		\hat u\\ \hat \theta
	\end{bmatrix} \in \mathscr X\times \mathscr S$ of the equation \eqref{MatrixEq} such that the norm $\norm{\begin{bmatrix}
			u_0-\hat u_0\\ \theta_0-\hat \theta_0
	\end{bmatrix}}_{L^p\cap L^2}$ is sufficiently small, we have 
	 \begin{eqnarray}\label{expstab}
	 	\norm{\begin{bmatrix}
	 			u-\hat u\\ \theta-\hat \theta
	 		\end{bmatrix}(t)}^\blacklozenge	\lesssim  e^{-\bf\Theta t} \norm{\begin{bmatrix}
	 		u_0-\hat u_0\\ \theta_0-\hat \theta_0
 		\end{bmatrix}}_{L^p\cap L^2}   \text{ for all } t>0,
	 \end{eqnarray} 
where $0<{\bf\Theta} <\beta$ (here $\beta$ is a positive constant and is given in Lemma \ref{estimates}).
\end{theorem} 
\begin{proof}
First, there exists a positive constant $\hat{\rho}$ such that $\norm{\begin{bmatrix}
		 \hat{u}\\ \hat{\theta}
	\end{bmatrix}}_{\mathscr{X}\times \mathscr{I}}<\hat{\rho}$.  
	
It is clear that
\begin{eqnarray}\label{ephi3}
	\begin{bmatrix}
			u-\hat  u\\ \theta-\hat \theta
		\end{bmatrix}(t)
\!\!\!&=&\!\!\! e^{-t\A}\begin{bmatrix}
		u_0-\hat u_0\\ \theta_0-\hat \theta_0
		\end{bmatrix}+ \mathscr B\left(\begin{bmatrix}
				u\\\theta
		\end{bmatrix},\begin{bmatrix}
				u\\\theta
		\end{bmatrix} \right)(t) - \mathscr B\left(\begin{bmatrix}
		\hat{u}\\\hat{\theta}
		\end{bmatrix},\begin{bmatrix}
		\hat{u}\\\hat{\theta}
		\end{bmatrix} \right)(t) + \cal H_h(\theta-\hat{\theta}) \cr
&\leq&\!\!\!e^{-t\A}\begin{bmatrix}
			u_0-\hat{u}_0\\ \theta_0-\hat{\theta}_0
		\end{bmatrix}+   \int_0^t e^{-(t-\tau)\A} \begin{bmatrix}
		\mathbb P[h(\theta-\hat{\theta}) \\ 0
		\end{bmatrix}(\tau)d\tau \cr&+&\!\!\!  \int_0^t e^{-(t-\tau)\A} \begin{bmatrix}
		\mathbb P([(u-\hat u) \cdot\nabla ] u+ [\hat u\cdot\nabla ](u-\hat u))\\ 
		 \;\,[(u-\hat u) \cdot \nabla ]\hat \theta+ [u\cdot \nabla](\theta-\hat \theta) 
		\end{bmatrix}(\tau)d\tau.  
	\end{eqnarray}
Hence,  we need to estimate $\norm{\begin{bmatrix}
		u-\hat u\\ \theta-\hat \theta
	\end{bmatrix}(t)}^\blacklozenge$  on the half time-line $t>0$. For this purpose, we frst have
\begin{eqnarray}\label{exp2}
&&\norm{\begin{bmatrix}
	u-\hat{u}\\ \theta-\hat{\theta}
	\end{bmatrix}(t)}_2 \leq  2Ce^{-\beta t}\norm{\begin{bmatrix}
	u_0-\hat{u}_0\\ \theta_0-\hat{\theta}_0
	\end{bmatrix}}_{2} + 2C\int_0^t 		  e^{-\beta(t-\tau)}\norm{\begin{bmatrix}
h(\theta-\hat{\theta}) \\ 0
\end{bmatrix}(\tau)}_{2}d\tau\cr
&&\hspace{3cm}+ 2C\int_0^t  e^{-\beta(t-\tau)}\norm{\begin{bmatrix}
		[(u-\hat u) \cdot\nabla ] u+ [\hat u\cdot\nabla ](u-\hat u)\\ 
	[(u-\hat u) \cdot \nabla ]\hat \theta+	[u\cdot \nabla](\theta-\hat \theta) 
	\end{bmatrix}(\tau)}_{2}d\tau\cr
&\leq& 2Ce^{-\beta t}\norm{\begin{bmatrix}
		u_0-\hat{u}_0\\ \theta_0-\hat{\theta}_0
\end{bmatrix}}_{2} + 2C\norm{h}_{\mathscr H}\int_0^t 		  e^{-\beta(t-\tau)}\norm{\begin{bmatrix}
0\\	\theta-\hat{\theta} 
\end{bmatrix}(\tau)}_{p}d\tau
\cr 
&&+ 2C\int_0^t  e^{-\beta(t-\tau)}\left[ \norm{
 (u-\hat{u})(\tau)}_p \norm{\begin{bmatrix}
\nabla u\\ \nabla \hat \theta
\end{bmatrix}(\tau)}_{\hat p} +\norm{
\begin{bmatrix}
	u \\ \hat{u}
	\end{bmatrix}}_p \norm{\begin{bmatrix}
\nabla (u-\hat u)\\ \nabla (\theta -\hat \theta)
\end{bmatrix}(\tau)}_{\hat p}  \right] d\tau  
\cr
&\leq& 2Ce^{-\beta t}\norm{\begin{bmatrix}
		u_0-\hat{u}_0\\ \theta_0-\hat{\theta}_0
\end{bmatrix}}_{2} + 2C\norm{h}_{\mathscr H}\int_0^t 		  e^{-\beta(t-\tau)}\norm{\begin{bmatrix}
		0\\	\theta-\hat{\theta} 
	\end{bmatrix}(\tau)}_{p}d\tau\cr 
&&+2C\left[ \norm{\begin{bmatrix}
		u\\
		\hat \theta
\end{bmatrix}}_{\mathscr X\times\mathscr S} +\norm{\begin{bmatrix}
		\hat u\\ 0
\end{bmatrix}}_{\mathscr X\times\mathscr S}\right]\int_0^t [h_d(\tau)]^{\frac{1}{p}-\frac{1}{\hat p}+\frac{1}{d}} e^{-\beta(t-\tau)} \cr
&&\hspace{1cm}\times \left[ \norm{
	(u-\hat{u})(\tau)}_p  + [h_d(\tau)]^{-(\frac{1}{p}-\frac{1}{\hat p}+\frac{1}{d})}\norm{\begin{bmatrix}
		\nabla (u-\hat u)\\ \nabla (\theta -\hat \theta)
	\end{bmatrix}(\tau)}_{\hat p}  \right] d\tau\cr
&\leq& 2Ce^{-\beta t}\norm{\begin{bmatrix}
		u_0-\hat{u}_0\\ \theta_0-\hat{\theta}_0
\end{bmatrix}}_{2} + 2C\norm{h}_{\mathscr H}\int_0^t 	[h_d(t-\tau)]^{\frac{1}{s}}	  e^{-\beta(t-\tau)}\norm{\begin{bmatrix}
		0\\	\theta-\hat{\theta} 
	\end{bmatrix}(\tau)}_{p}d\tau
\cr &&+2C(\rho+\hat{\rho}) \int_0^t [h_d(\tau)]^{\frac{1}{p}-\frac{1}{\hat p}+\frac{1}{d}} e^{-\beta(t-\tau)} 
\norm{\begin{bmatrix}
		u-\hat u\\ \theta -\hat \theta
	\end{bmatrix}(\tau)}^\blacklozenge d\tau,
\end{eqnarray}
where, we used $[h_d(t-\tau)]^{\frac{1}{s}}>1,\; \forall t>\tau$.  
   
From the condition $2\leq d<p<\hat p< s$, we have $1<[h_d(\tau)]^{\frac{1}{p}-\frac{1}{\hat p }+\frac{1}{d}} < [h_d(\tau)]^{\frac{1}{p}-\frac{1}{s}+\frac{1}{d}}$ for all $\tau>0.$  Then, we can estimate   
\begin{eqnarray}\label{expp}
&&\norm{\begin{bmatrix}
   			u-\hat{u}\\ \theta-\hat{\theta}
   		\end{bmatrix}(t)}_p \leq  Ce^{-\beta t}\norm{\begin{bmatrix}
   			u_0-\hat{u}_0\\ \theta_0-\hat{\theta}_0
   	\end{bmatrix}}_{L^p\cap L^2}  + C\int_0^t [h_d(t-\tau)]^{\frac{1}{r}-\frac{1}{ p}} e^{-\beta(t-\tau)}\norm{\begin{bmatrix}
   	h(\theta-\hat{\theta}) \\ 0
   	\end{bmatrix}(\tau)}_{L^r\cap L^2}d\tau\cr
&&\hspace{3cm}+ C\int_0^t [h_d(t-\tau)]^{\frac{1}{r}-\frac{1}{ p}} e^{-\beta(t-\tau)}\norm{\begin{bmatrix}
   			[(u-\hat u) \cdot\nabla ] u+ [\hat u\cdot\nabla](u-\hat u)\\ 
   			[(u-\hat u) \cdot \nabla ]\hat \theta+	[u\cdot \nabla](\theta-\hat \theta) 
   		\end{bmatrix}(\tau)}_{L^r\cap L^2} d\tau \cr
&\leq& Ce^{-\beta t}\norm{\begin{bmatrix}
   			u_0-\hat{u}_0\\ \theta_0-\hat{\theta}_0
   	\end{bmatrix}}_{L^p\cap L^2} + C\norm{h}_{\mathscr H}\int_0^t 	[h_d(t-\tau)]^{\frac{1}{s}}  e^{-\beta(t-\tau)}\norm{\begin{bmatrix}
   			0\\	\theta-\hat{\theta} 
   		\end{bmatrix}(\tau)}_{p}d\tau
   	\cr &&+ C\int_0^t [h_d(t-\tau)]^{\frac{1}{s}} e^{-\beta(t-\tau)}\left[\norm{(u-\hat{u})(\tau)}_p \norm{\begin{bmatrix}
   			\nabla u\\ \nabla \hat \theta
   		\end{bmatrix}(\tau)}_{L^{\hat p}\cap L^{s}} +\norm{
   		\begin{bmatrix}
   			u \\ \hat{u}
   	\end{bmatrix}}_p \norm{\begin{bmatrix}
   			\nabla (u-\hat u)\\ \nabla (\theta -\hat \theta)
   		\end{bmatrix}(\tau)}_{L^{\hat p}\cap L^{s}} \right] d\tau  \cr
&\leq& Ce^{-\beta t}\norm{\begin{bmatrix}
   			u_0-\hat{u}_0\\ \theta_0-\hat{\theta}_0
   	\end{bmatrix}}_{L^p\cap L^2} + C\norm{h}_{\mathscr H}\int_0^t 	[h_d(t-\tau)]^{\frac{1}{s}}	  e^{-\beta(t-\tau)}\norm{\begin{bmatrix}
   			0\\	\theta-\hat{\theta} 
   		\end{bmatrix}(\tau)}_{p}d\tau
   	\cr 
&&+C(\rho+\hat\rho)\int_0^t[h_d(t-\tau)]^{\frac{1}{s}} [h_d(\tau)]^{\frac{1}{p}-\frac{1}{\hat p}+\frac{1}{d}}   e^{-\beta(t-\tau)}\left[ \norm{
   		(u-\hat{u})(\tau)}_p   + [h_d(\tau)]^{-(\frac{1}{p}-\frac{1}{\hat p}+\frac{1}{d})}\norm{\begin{bmatrix}
   	\nabla (u-\hat u)\\ \nabla (\theta -\hat \theta)
   		\end{bmatrix}(\tau)}_{\hat p}  \right] d\tau
   	\cr &&+C(\rho+\hat\rho)\int_0^t [h_d(t-\tau)]^{\frac{1}{s}} [h_d(\tau)]^{\frac{1}{p}-\frac{1}{s}+\frac{1}{d}} e^{-\beta(t-\tau)} \left[ \norm{
   		(u-\hat{u})(\tau)}_p  + [h_d(\tau)]^{-(\frac{1}{p}-\frac{1}{s}+\frac{1}{d})}\norm{\begin{bmatrix}
   			\nabla (u-\hat u)\\ \nabla (\theta -\hat \theta)
   		\end{bmatrix}(\tau)}_{s} \right] d\tau
\cr
&\leq& Ce^{-\beta t}\norm{\begin{bmatrix}
		u_0-\hat{u}_0\\ \theta_0-\hat{\theta}_0
\end{bmatrix}}_{L^p\cap L^2} + C\norm{h}_{\mathscr H}\int_0^t 	[h_d(t-\tau)]^{\frac{1}{s}}	  e^{-\beta(t-\tau)}\norm{\begin{bmatrix}
		0\\	\theta-\hat{\theta} 
	\end{bmatrix}(\tau)}_{p}d\tau
\cr &&+2C(\rho+\hat\rho)\int_0^t [h_d(t-\tau)]^{\frac{1}{s}} [h_d(\tau)]^{\frac{1}{p}-\frac{1}{s}+\frac{1}{d}} e^{-\beta(t-\tau)}  \norm{\begin{bmatrix}
		u-\hat u\\ \theta -\hat \theta
	\end{bmatrix}(\tau)}^\blacklozenge d\tau.
   \end{eqnarray}
  
Moreover, we also have the following gradient estimate
\begin{eqnarray}\label{expp^}
&&[h_d(t)]^{-(\frac{1}{p}-\frac{1}{\hat p}+\frac{1}{d})}\norm{\begin{bmatrix}
			\nabla (u-\hat u)\\ \nabla (\theta -\hat \theta)
		\end{bmatrix}(t)}_{\hat p}\leq  Ce^{-\beta t}\norm{\begin{bmatrix}
	u_0-\hat{u}_0\\ \theta_0-\hat{\theta}_0
	\end{bmatrix}}_{L^p\cap L^2} 
\cr&& + C\int_0^t[h_d(t)]^{-(\frac{1}{p}-\frac{1}{\hat p}+\frac{1}{d})} [h_d(t-\tau)]^{\frac{1}{r}-\frac{1}{\hat  p}+\frac{1}{d}} e^{-\beta(t-\tau)}\norm{\begin{bmatrix}
			h(\theta-\hat{\theta}) \\ 0
		\end{bmatrix}(\tau)}_{L^r\cap L^2}d\tau\cr
&&+ C\int_0^t[h_d(t)]^{-(\frac{1}{p}-\frac{1}{\hat p}+\frac{1}{d})} [h_d(t-\tau)]^{\frac{1}{r}-\frac{1}{\hat p}+\frac{1}{d}} e^{-\beta(t-\tau)} \norm{\begin{bmatrix}
			[(u-\hat u) \cdot\nabla ] u+ [\hat u\cdot\nabla](u-\hat u)\\ 
			[(u-\hat u) \cdot \nabla ]\hat \theta+	[u\cdot \nabla](\theta-\hat \theta) 
		\end{bmatrix}(\tau)}_{L^r\cap L^2} d\tau \cr
&\leq& Ce^{-\beta t}\norm{\begin{bmatrix}
			u_0-\hat{u}_0\\ \theta_0-\hat{\theta}_0
	\end{bmatrix}}_{L^p\cap L^2} + C\norm{h}_{\mathscr H}\int_0^t 	[h_d(t)]^{-(\frac{1}{p}-\frac{1}{\hat p}+\frac{1}{d})} [h_d(t-\tau)]^{\frac{1}{r}-\frac{1}{\hat  p}+\frac{1}{d}}  e^{-\beta(t-\tau)}\norm{\begin{bmatrix}
			0\\	\theta-\hat{\theta} 
		\end{bmatrix}(\tau)}_{p}d\tau
	\cr 
&&+ C\int_0^t[h_d(t)]^{-(\frac{1}{p}-\frac{1}{\hat p}+\frac{1}{d})} [h_d(t-\tau)]^{\frac{1}{r}-\frac{1}{\hat  p}+\frac{1}{d}} e^{-\beta(t-\tau)}\cr&&\hspace{1cm}\times\left[\norm{
		(u-\hat{u})(\tau)}_p \norm{\begin{bmatrix}
			\nabla u\\ \nabla \hat \theta
		\end{bmatrix}(\tau)}_{L^{\hat p}\cap L^{s}} +\norm{
		\begin{bmatrix}
			u \\ \hat{u}
	\end{bmatrix}}_p \norm{\begin{bmatrix}
			\nabla (u-\hat u)\\ \nabla (\theta -\hat \theta)
		\end{bmatrix}(\tau)}_{L^{\hat p}\cap L^{s}} \right] d\tau 
\cr
&\leq& Ce^{-\beta t}\norm{\begin{bmatrix}
			u_0-\hat{u}_0\\ \theta_0-\hat{\theta}_0
	\end{bmatrix}}_{L^p\cap L^2} + C\norm{h}_{\mathscr H}\int_0^t [h_d(t)]^{-(\frac{1}{p}-\frac{1}{\hat p}+\frac{1}{d})} [h_d(t-\tau)]^{\frac{1}{r}-\frac{1}{\hat  p}+\frac{1}{d}}	  e^{-\beta(t-\tau)}\norm{\begin{bmatrix}
			0\\	\theta-\hat{\theta} 
		\end{bmatrix}(\tau)}_{p}d\tau\cr 
&&+C(\rho+\hat{\rho})\int_0^t[h_d(t)]^{-(\frac{1}{p}-\frac{1}{\hat p}+\frac{1}{d})} [h_d(t-\tau)]^{\frac{1}{r}-\frac{1}{\hat  p}+\frac{1}{d}} [h_d(\tau)]^{\frac{1}{p}-\frac{1}{\hat p}+\frac{1}{d}}   e^{-\beta(t-\tau)} \cr
&&\hspace{1cm}\times \left[ \norm{
		(u-\hat{u})(\tau)}_p   + [h_d(\tau)]^{-(\frac{1}{p}-\frac{1}{\hat p}+\frac{1}{d})}\norm{\begin{bmatrix}
			\nabla (u-\hat u)\\ \nabla (\theta -\hat \theta)
		\end{bmatrix}(\tau)}_{\hat p}  \right] d\tau\cr 
&&+C(\rho+\hat\rho)\int_0^t [h_d(t)]^{-(\frac{1}{p}-\frac{1}{\hat p}+\frac{1}{d})} [h_d(t-\tau)]^{\frac{1}{r}-\frac{1}{\hat  p}+\frac{1}{d}} [h_d(\tau)]^{\frac{1}{p}-\frac{1}{s}+\frac{1}{d}} e^{-\beta(t-\tau)} \cr&&\hspace{1cm}\times \left[ \norm{
		(u-\hat{u})(\tau)}_p  + [h_d(\tau)]^{-(\frac{1}{p}-\frac{1}{s}+\frac{1}{d})}\norm{\begin{bmatrix}
			\nabla (u-\hat u)\\ \nabla (\theta -\hat \theta)
		\end{bmatrix}(\tau)}_{s} \right] d\tau\cr
&\leq& Ce^{-\beta t}\norm{\begin{bmatrix}
		u_0-\hat{u}_0\\ \theta_0-\hat{\theta}_0
\end{bmatrix}}_{L^p\cap L^2} + C\norm{h}_{\mathscr H}\int_0^t [h_d(t)]^{-(\frac{1}{p}-\frac{1}{\hat p}+\frac{1}{d})} [h_d(t-\tau)]^{\frac{1}{r}-\frac{1}{\hat  p}+\frac{1}{d}}	  e^{-\beta(t-\tau)}\norm{\begin{bmatrix}
		0\\	\theta-\hat{\theta} 
	\end{bmatrix}(\tau)}_{p}d\tau\cr
&&+2C(\rho+\hat\rho)\int_0^t [h_d(t)]^{-(\frac{1}{p}-\frac{1}{\hat p}+\frac{1}{d})} [h_d(t-\tau)]^{\frac{1}{r}-\frac{1}{\hat  p}+\frac{1}{d}} [h_d(\tau)]^{\frac{1}{p}-\frac{1}{s}+\frac{1}{d}} e^{-\beta(t-\tau)}  \norm{\begin{bmatrix}
		u-\hat u\\  \theta -\hat \theta
	\end{bmatrix}(\tau)}^\blacklozenge d\tau.\cr&&
\end{eqnarray}   
Finally, we have   
\begin{eqnarray}\label{expsigma}
&&[h_d(t)]^{-(\frac{1}{p}-\frac{1}{s}+\frac{1}{d})}\norm{\begin{bmatrix}
			\nabla (u-\hat u)\\ \nabla (\theta -\hat \theta)
		\end{bmatrix}(t)}_{s} \leq  Ce^{-\beta t}\norm{\begin{bmatrix}
			u_0-\hat{u}_0\\ \theta_0-\hat{\theta}_0
	\end{bmatrix}}_{L^p\cap L^2} \cr
&& + C\int_0^t[h_d(t)]^{-(\frac{1}{p}-\frac{1}{s}+\frac{1}{d})} [h_d(t-\tau)]^{\frac{1}{r}-\frac{1}{s}+\frac{1}{d}} e^{-\beta(t-\tau)}\norm{\begin{bmatrix}
			h(\theta-\hat{\theta}) \\ 0
		\end{bmatrix}(\tau)}_{L^r\cap L^2}d\tau\cr
	&&+C\int_0^t[h_d(t)]^{-(\frac{1}{p}-\frac{1}{s}+\frac{1}{d})} [h_d(t-\tau)]^{\frac{1}{r}-\frac{1}{s}+\frac{1}{d}} e^{-\beta(t-\tau)} \norm{\begin{bmatrix}
			[(u-\hat u) \cdot\nabla ] u+ [\hat u\cdot\nabla](u-\hat u)\\ 
			[(u-\hat u) \cdot \nabla ]\hat \theta+	[u\cdot \nabla](\theta-\hat \theta) 
		\end{bmatrix}(\tau)}_{L^r\cap L^2} d\tau
\cr
&\leq& Ce^{-\beta t}\norm{\begin{bmatrix}
			u_0-\hat{u}_0\\ \theta_0-\hat{\theta}_0
	\end{bmatrix}}_{L^p\cap L^2} 
 + C\norm{h}_{\mathscr H}\int_0^t 	[h_d(t)]^{-(\frac{1}{p}-\frac{1}{s}+\frac{1}{d})} [h_d(t-\tau)]^{\frac{1}{r}-\frac{1}{s}+\frac{1}{d}}  e^{-\beta(t-\tau)}\norm{\begin{bmatrix}
			0\\	\theta-\hat{\theta} 
		\end{bmatrix}(\tau)}_{p}d\tau
	\cr &&+ C\int_0^t[h_d(t)]^{-(\frac{1}{p}-\frac{1}{s}+\frac{1}{d})} [h_d(t-\tau)]^{\frac{1}{r}-\frac{1}{s}+\frac{1}{d}} e^{-\beta(t-\tau)}\cr&&\hspace{1cm}\times\left[\norm{
		(u-\hat{u})(\tau)}_p \norm{\begin{bmatrix}
			\nabla u\\ \nabla \hat \theta
		\end{bmatrix}(\tau)}_{L^{\hat p}\cap L^{s}} +\norm{
		\begin{bmatrix}
			u \\ \hat{u}
	\end{bmatrix}}_p \norm{\begin{bmatrix}
			\nabla (u-\hat u)\\ \nabla (\theta -\hat \theta)
		\end{bmatrix}(\tau)}_{L^{\hat p}\cap L^{s}} \right] d\tau 
\cr
&\leq& Ce^{-\beta t}\norm{\begin{bmatrix}
			u_0-\hat{u}_0\\ \theta_0-\hat{\theta}_0
	\end{bmatrix}}_{L^p\cap L^2} + C\norm{h}_{\mathscr H}\int_0^t [h_d(t)]^{-(\frac{1}{p}-\frac{1}{s}+\frac{1}{d})} [h_d(t-\tau)]^{\frac{1}{r}-\frac{1}{s}+\frac{1}{d}}	  e^{-\beta(t-\tau)}\norm{\begin{bmatrix}
			0\\	\theta-\hat{\theta} 
		\end{bmatrix}(\tau)}_{p}d\tau\cr &&+C(\rho+\hat\rho)\int_0^t[h_d(t)]^{-(\frac{1}{p}-\frac{1}{s}+\frac{1}{d})} [h_d(t-\tau)]^{\frac{1}{r}-\frac{1}{s}+\frac{1}{d}} [h_d(\tau)]^{\frac{1}{p}-\frac{1}{\hat p}+\frac{1}{d}}   e^{-\beta(t-\tau)} \cr&&\hspace{1cm}\times \left[ \norm{
		(u-\hat{u})(\tau)}_p   + [h_d(\tau)]^{-(\frac{1}{p}-\frac{1}{\hat p}+\frac{1}{d})}\norm{\begin{bmatrix}
			\nabla (u-\hat u)\\ \nabla (\theta -\hat \theta)
		\end{bmatrix}(\tau)}_{\hat p}  \right] d\tau	\cr 
&&+C(\rho+\hat\rho)\int_0^t [h_d(t)]^{-(\frac{1}{p}-\frac{1}{s}+\frac{1}{d})} [h_d(t-\tau)]^{\frac{1}{r}-\frac{1}{s}+\frac{1}{d}} [h_d(\tau)]^{\frac{1}{p}-\frac{1}{s}+\frac{1}{d}} e^{-\beta(t-\tau)} \cr&&\hspace{1cm}\times \left[ \norm{
		(u-\hat{u})(\tau)}_p  + [h_d(\tau)]^{-(\frac{1}{p}-\frac{1}{s}+\frac{1}{d})}\norm{\begin{bmatrix}
			\nabla (u-\hat u)\\ \nabla (\theta -\hat \theta)
		\end{bmatrix}(\tau)}_{s} \right] d\tau
\cr
&\leq& Ce^{-\beta t}\norm{\begin{bmatrix}
		u_0-\hat{u}_0\\ \theta_0-\hat{\theta}_0
\end{bmatrix}}_{L^p\cap L^2} + C\norm{h}_{\mathscr H}\int_0^t [h_d(t)]^{-(\frac{1}{p}-\frac{1}{s}+\frac{1}{d})} [h_d(t-\tau)]^{\frac{1}{r}-\frac{1}{s}+\frac{1}{d}}	  e^{-\beta(t-\tau)}\norm{\begin{bmatrix}
		0\\	\theta-\hat{\theta} 
	\end{bmatrix}(\tau)}_{p}d\tau
\cr &&+2C(\rho+\hat\rho)\int_0^t [h_d(t)]^{-(\frac{1}{p}-\frac{1}{s}+\frac{1}{d})} [h_d(t-\tau)]^{\frac{1}{r}-\frac{1}{s}+\frac{1}{d}} [h_d(\tau)]^{\frac{1}{p}-\frac{1}{s}+\frac{1}{d}} e^{-\beta(t-\tau)}  \norm{\begin{bmatrix}
		u-\hat u\\ \theta -\hat \theta
	\end{bmatrix}(\tau)}^\blacklozenge d\tau.
\end{eqnarray}    
   
Combining \eqref{exp2}, \eqref{expp}, \eqref{expp^} and \eqref{expsigma}, we deduce that 
 \begin{eqnarray*}\label{expsum1}
\norm{\begin{bmatrix}
		u-\hat u\\ \theta -\hat \theta
		\end{bmatrix}(t)}^\blacklozenge
&\leq&  5Ce^{-\beta t}\norm{\begin{bmatrix}
			u_0-\hat{u}_0\\ \theta_0-\hat{\theta}_0
	\end{bmatrix}}_{L^p\cap L^2} 
+2C(\rho+\hat\rho) \int_0^t [h_d(\tau)]^{\frac{1}{p}-\frac{1}{\hat p}+\frac{1}{d}} e^{-\beta(t-\tau)} 
\norm{\begin{bmatrix}
		u-\hat u\\ \theta -\hat \theta
	\end{bmatrix}(\tau)}^\blacklozenge d\tau
\cr
&& +3 C\norm{h}_{\mathscr H}\int_0^t 	[h_d(t-\tau)]^{\frac{1}{s}}	  e^{-\beta(t-\tau)}\norm{\begin{bmatrix}
		0\\	\theta-\hat{\theta} 
	\end{bmatrix}(\tau)}_{p}d\tau\cr 
&&+2C(\rho+\hat\rho)\int_0^t [h_d(t-\tau)]^{\frac{1}{s}} [h_d(\tau)]^{\frac{1}{p}-\frac{1}{s}+\frac{1}{d}} e^{-\beta(t-\tau)}  \norm{\begin{bmatrix}
		u-\hat u\\ \theta -\hat \theta
	\end{bmatrix}(\tau)}^\blacklozenge d\tau
\cr
&& 
+ C\norm{h}_{\mathscr H}\int_0^t [h_d(t)]^{-(\frac{1}{p}-\frac{1}{\hat p}+\frac{1}{d})} [h_d(t-\tau)]^{\frac{1}{r}-\frac{1}{\hat  p}+\frac{1}{d}}	  e^{-\beta(t-\tau)}\norm{\begin{bmatrix}
		0\\	\theta-\hat{\theta} 
	\end{bmatrix}(\tau)}_{p}d\tau\cr 
&&+2C(\rho+\hat\rho)\int_0^t [h_d(t)]^{-(\frac{1}{p}-\frac{1}{\hat p}+\frac{1}{d})} [h_d(t-\tau)]^{\frac{1}{r}-\frac{1}{\hat  p}+\frac{1}{d}} [h_d(\tau)]^{\frac{1}{p}-\frac{1}{s}+\frac{1}{d}} e^{-\beta(t-\tau)}  \norm{\begin{bmatrix}
		u-\hat u\\  \theta -\hat \theta
	\end{bmatrix}(\tau)}^\blacklozenge d\tau \cr&&+ C\norm{h}_{\mathscr H}\int_0^t [h_d(t)]^{-(\frac{1}{p}-\frac{1}{s}+\frac{1}{d})} [h_d(t-\tau)]^{\frac{1}{r}-\frac{1}{s}+\frac{1}{d}}	  e^{-\beta(t-\tau)}\norm{\begin{bmatrix}
			0\\	\theta-\hat{\theta} 
		\end{bmatrix}(\tau)}_{p}d\tau\cr &&+2C(\rho+\hat\rho)\int_0^t [h_d(t)]^{-(\frac{1}{p}-\frac{1}{s}+\frac{1}{d})} [h_d(t-\tau)]^{\frac{1}{r}-\frac{1}{s}+\frac{1}{d}} [h_d(\tau)]^{\frac{1}{p}-\frac{1}{s}+\frac{1}{d}} e^{-\beta(t-\tau)}  \norm{\begin{bmatrix}
			u-\hat u\\ \theta -\hat \theta
		\end{bmatrix}(\tau)}^\blacklozenge d\tau.
\end{eqnarray*}

From the fact that  $1<[h_d(t-\tau)]^{\frac{1}{s}}$, and $1<[h_d(\tau)]^{\frac{1}{p}-\frac{1}{\hat p }+\frac{1}{d}} < [h_d(\tau)]^{\frac{1}{p}-\frac{1}{s}+\frac{1}{d}}$ for all $t>\tau>0$, where $2\leq d<p<\hat p< s$, we have
\begin{eqnarray*}\label{expsum2}
&&\norm{\begin{bmatrix}
			u-\hat u\\ \theta -\hat \theta
		\end{bmatrix}(t)}^\blacklozenge
\leq  5Ce^{-\beta t}\norm{\begin{bmatrix}
	u_0-\hat{u}_0\\ \theta_0-\hat{\theta}_0
	\end{bmatrix}}_{L^p\cap L^2} \cr 
&&+C\left(4(\rho+\hat\rho)+3\norm{h}_{\mathscr H}\right) \int_0^t [h_d(t-\tau)]^{\frac{1}{s}} [h_d(\tau)]^{\frac{1}{p}-\frac{1}{s}+\frac{1}{d}} e^{-\beta(t-\tau)}  \norm{\begin{bmatrix}
			u-\hat u\\ \theta -\hat \theta
		\end{bmatrix}(\tau)}^\blacklozenge d\tau
\cr &&	+C\left(2(\rho+\hat\rho)+\norm{h}_{\mathscr H}\right) \int_0^t [h_d(t)]^{-(\frac{1}{p}-\frac{1}{\hat p}+\frac{1}{d})} [h_d(t-\tau)]^{\frac{1}{r}-\frac{1}{\hat  p}+\frac{1}{d}} [h_d(\tau)]^{\frac{1}{p}-\frac{1}{s}+\frac{1}{d}} e^{-\beta(t-\tau)}  \norm{\begin{bmatrix}
	u-\hat u\\  \theta -\hat \theta
	\end{bmatrix}(\tau)}^\blacklozenge d\tau 
\cr 
&&+C\left(2(\rho+\hat\rho)+\norm{h}_{\mathscr H}\right) \int_0^t [h_d(t)]^{-(\frac{1}{p}-\frac{1}{s}+\frac{1}{d})} [h_d(t-\tau)]^{\frac{1}{r}-\frac{1}{s}+\frac{1}{d}} [h_d(\tau)]^{\frac{1}{p}-\frac{1}{s}+\frac{1}{d}} e^{-\beta(t-\tau)}  \norm{\begin{bmatrix}
			u-\hat u\\ \theta -\hat \theta
		\end{bmatrix}(\tau)}^\blacklozenge d\tau.
\end{eqnarray*}

Setting $z(\tau)=e^{\bf\Theta \tau} \norm{\begin{bmatrix}
		u-\tilde{u}\\ \theta-\tilde{\theta}
	\end{bmatrix}(\tau)}^\blacklozenge $ for ${\bf \Theta}< \beta$. From above estimations, we have
\begin{eqnarray}\label{expsum3}
z(t)&\leq&  5Ce^{-(\beta-{\bf \Theta}) t}\norm{\begin{bmatrix}
			u_0-\hat{u}_0\\ \theta_0-\hat{\theta}_0
	\end{bmatrix}}_{L^p\cap L^2} \cr &&+C\left(4(\rho+\hat\rho)+3\norm{h}_{\mathscr H}\right) \int_0^t [h_d(t-\tau)]^{\frac{1}{s}} [h_d(\tau)]^{\frac{1}{p}-\frac{1}{s}+\frac{1}{d}} e^{-(\beta-{\bf \Theta})(t-\tau)}  z(\tau)d\tau\cr 
&&+C\left(2(\rho+\hat\rho)+\norm{h}_{\mathscr H}\right) \int_0^t [h_d(t)]^{-(\frac{1}{p}-\frac{1}{\hat p}+\frac{1}{d})} [h_d(t-\tau)]^{\frac{1}{r}-\frac{1}{\hat  p}+\frac{1}{d}} [h_d(\tau)]^{\frac{1}{p}-\frac{1}{s}+\frac{1}{d}} e^{-(\beta-{\bf \Theta})(t-\tau)}  z(\tau)d\tau \cr &&+C\left(2(\rho+\hat\rho)+\norm{h}_{\mathscr H}\right) \int_0^t [h_d(t)]^{-(\frac{1}{p}-\frac{1}{s}+\frac{1}{d})} [h_d(t-\tau)]^{\frac{1}{r}-\frac{1}{s}+\frac{1}{d}} [h_d(\tau)]^{\frac{1}{p}-\frac{1}{s}+\frac{1}{d}} e^{-(\beta-{\bf \Theta})(t-\tau)}  z(\tau)d\tau.
\end{eqnarray}

Similar to the boundedness of integrals by constants $K_p, K_{\hat p}$ and $K_s$ in the proof of Lemma \ref{BiE}, we can also have the following convergence
$$ \int_0^t [h_d(t-\tau)]^{\frac{1}{s}} [h_d(\tau)]^{\frac{1}{p}-\frac{1}{s}+\frac{1}{d}} e^{-(\beta-{\bf \Theta})(t-\tau)}  z(\tau)d\tau <\infty,$$
$$ \int_0^t [h_d(t)]^{-(\frac{1}{p}-\frac{1}{\hat p}+\frac{1}{d})} [h_d(t-\tau)]^{\frac{1}{r}-\frac{1}{\hat  p}+\frac{1}{d}} [h_d(\tau)]^{\frac{1}{p}-\frac{1}{s}+\frac{1}{d}} e^{-(\beta-{\bf \Theta})(t-\tau)}  z(\tau)d\tau <\infty,$$
and
$$ \int_0^t [h_d(t)]^{-(\frac{1}{p}-\frac{1}{s}+\frac{1}{d})} [h_d(t-\tau)]^{\frac{1}{r}-\frac{1}{s}+\frac{1}{d}} [h_d(\tau)]^{\frac{1}{p}-\frac{1}{s}+\frac{1}{d}} e^{-(\beta-{\bf \Theta})(t-\tau)}  z(\tau)d\tau <\infty.$$
Therefore, we use Gronwall's inequality for \eqref{expsum3} to get 
	$$|z(t)|\leq \widehat{C} \norm{\begin{bmatrix}
			u_0-\hat{u}_0\\ \theta_0-\hat{\theta}_0
	\end{bmatrix}}_{L^p\cap L^2} \text{ for all }t>0.$$
This implies the desired stability \eqref{expstab}. 
\end{proof}

\section{Generalized gravitational fields and some boundedness of integrals}
\subsection{Generalized gravitational fields on non-Euclidean frameworks}\label{appen}
In a recent work \cite{XuanTrung2022}, we have consider the generalized gravitional fields on real hyperbolic spaces satisfying Assumption \ref{Assum}. Now, we discuss about the existence of these fields on the generalized non-compact Riemannian manifolds.
In particular, we consider a general static spherically symmetric Riemannian manifold $({\bf M},g)$ with dimension $d$, where the metric $g$ is given by 
$$g = dr^2 + f(r)d\omega^2.$$
Here, $d\omega^2$ denotes the standard Euclidean metric on $(d-1)$ dimension sphere $\mathbb{S}^{d-1}$ and $f(r)$ is a smooth positive function of $r$. The manifold $({\bf M},g)$ is a case of Cartan-Hadamard manifolds and it satisfies conditions $(H_1)-(H_4)$ in Assumption \ref{Ass} if we restrict some conditions on function $f(r)$ to guarantee the boundedness and negative properties of curvatures. 

Now, we show that the gravitational field associated with the metric $g$ can be extended to a function satisfying Assumption \ref{Assum}. First, the gravitational potential $\Phi(r)$ associated with metric $g$ is a solution of Poisson equation
\begin{equation*}
\nabla^2\Phi(r) = 0.
\end{equation*}
The solution is (see \cite[page 151]{Go} and also \cite[Section 5]{Barrow2020}):
\begin{equation*}
\Phi(r) = G(x,O) = \int_r^\infty \frac{d\rho}{f(\rho)},
\end{equation*}
where $G$ is the Green function, $O$ is the origin point of manifold and $\rho=\mathrm{dist}(x,O)$. 

If ${\bf M}=\mathbb{R}^d$ (where $d\geqslant 3$), then $f(r)=r^2$ and by a straightforward calculation we have (see \cite{Go}):
$$\Phi(r) = \frac{((d-2)\omega_d)^{-1}}{|r|^{d-2}},$$ 
where $\omega_d$ is the area of the unit sphere in $\mathbb{R}^d$. Then, the gravitational field is proportional to 
$$h(x) \simeq \nabla\Phi(r) = \frac{d}{dr}\frac{((d-2)\omega_d)^{-1}}{|r|^{d-2}}= - \frac{1}{\omega_d}\frac{r}{|r|^d}$$
by multiplying a positive constant depending on gravitational constant and mass.
Therefore, the function $h$ belongs to $L^\infty(T\Omega)\cap L^{\frac{d}{2},\infty}(T\Omega)$ (see \cite{Hi1997}), where $\Omega=\mathbb{R}^d-B(O,\varepsilon)$ is an exterior domain of $B(O,\varepsilon)$ which is a ball with origin at $O=(0,...,0)$ and radius $\varepsilon>0$.

If ${\bf M}=\mathbb{H}^d$ (where $d\geqslant 3$), then $f(r)=\sinh r$ and we have (see \cite{Barrow2020}):
$$\Phi(r) = \coth^{d-2} r$$
This leads to
$$h(x) \simeq \frac{d\Phi(r)}{d r} = 
-\frac{(d-2)\coth^{d-3}r}{\sinh^2 r}.$$
Observe that, for $d\geqslant 3$, we have the following equivalence
\begin{equation}\label{gra1}
{h}(x) \simeq
\begin{cases}
-e^{-2r}, \hbox{  as  } r\to \infty,\cr
-r^{-2}, \hbox{  as } r\to 0.
\end{cases}
\end{equation}
Now, we verify ${h}\in L^\infty(\Gamma(T\Omega)) \cap L^{\frac{d}{2},\infty}(\Gamma(T\Omega))$, for $h$ given by \eqref{gra1}, where $\Omega = {\bf M}-\mathbb{B}(O,\varepsilon)$ is an exterior of  $\mathbb{B}(O,\varepsilon)$ which is a geodesic ball in $\mathbb{H}^d$ centered at the origin $O= (1,0,0...0)$ with geodesic radius $\varepsilon$. Clearly, on the exterior domain $\Omega$, the condition $\tilde{h} \in L^{\infty}(\Gamma(T\Omega))$ is valid. Moreover, the norm of interpolation space defined on hyperbolic manifold ${\bf M}$ is given by (see notions in the proof of Corollary 3.3 in \cite{An2009}):
$$\norm{f}_{L^{q,\infty}} = \sup_{0<r<1}r^{\frac{d}{q}}|f(r)| + \sup_{r\geqslant 1} e^{\frac{d-1}{q}r}|f(r)|.$$
Using equivalence \eqref{gra1}, we can show that the gravitational field $\tilde{h}$ on ${\bf M}$ (with dimension $d\geqslant 3$) satisfies $\norm{{h}}_{L^{\frac{d}{2},\infty}}<+\infty$, hence ${h}$ belongs to $L^{\frac{d}{2},\infty}(\Gamma(T\Omega))$.

Therefore, if we consider other manifolds with $f(r)$ decays fastly as $r^\alpha$ (where $\alpha\geqslant 2$) and $e^{\beta r}$ (where $\beta \geqslant 1$), then we can also get that $h$ belongs to $L^\infty(T\Omega)\cap L^{\frac{d}{2},\infty}(T\Omega)$, where $\Omega$ is an exterior domain of a geodesic ball in ${\bf M}$ centered at the origin $O$ of ${\bf M}$.
With the function $h(x)$ on hand, we can construct a generalized gravitational field as $\widetilde{h}(x)=\gamma(x,t)h(x)$, where $\gamma:{\bf M}\times \mathbb{R}_+\to \mathbb{R}$, is a bounded and continuous function and has support outside the geodesic ball $\mathbb{B}(O,\varepsilon)$. Since ${h} \in L^\infty(\Gamma(T\Omega)) \cap L^{\frac{d}{2},\infty}(\Gamma(T\Omega))$ and the property of $\alpha(x,t)$, the generalized gravitational field $\widetilde{h}(x,t)=\alpha(x,t){h}(x)$ satisfies Assumption \ref{Assum}.

\subsection{Some boundedness of integrals}\label{convergence}
We give the details of the boundedness of $G_p(t)$ and $G_s(t)$ used in the proof of Lemma \ref{BiE}. Recall that
$$G_p(t) = C\int_0^t  [h_d(t-\tau)]^{\frac{1}{s}}\left[ [h_d(\tau)]^{\frac{1}{p}-\frac{1}{s}+\frac{1}{d}} +[h_d(\tau)]^{\frac{1}{p}-\frac{1}{\hat p}+\frac{1}{d}} \right]e^{-\beta(t-\tau)}  d\tau $$
Since $2\leq d<p<\hat p< s$, we have $\frac{1}{p}-\frac{1}{s}+\frac{1}{d} >\frac{1}{p}-\frac{1}{\hat p}+\frac{1}{d}>0$ and then $ [h_d(\tau)]^{\frac{1}{p}-\frac{1}{s}+\frac{1}{d}} >[h_d(\tau)]^{\frac{1}{p}-\frac{1}{\hat p}+\frac{1}{d}} $ for all $\tau>0$. Hence, we can estimate
 \begin{eqnarray*}
 G_p(t)\leq 2C\int_0^t  [h_d(t-\tau)]^{\frac{1}{s}} [h_d(\tau)]^{\frac{1}{p}-\frac{1}{s}+\frac{1}{d}} e^{-\beta(t-\tau)}  d\tau.
 \end{eqnarray*}
For the case $0<t\leq 1$, we have clearly that	 
 \begin{eqnarray}\label{Gp1}
	G_p(t)&\leq& 2C\int_0^t  (t-\tau)^{-\frac{d}{2s}} \tau^{-\frac{d}{2}
	(\frac{1}{p}-\frac{1}{s}+\frac{1}{d})}  d\tau \leq 2C\int_0^t  \left(1 -\frac{\tau}{t}\right) ^{-\frac{d}{2s}} \left( \frac{\tau}{t}\right) ^{-\frac{d}{2}
	(\frac{1}{p}-\frac{1}{s}+\frac{1}{d})}  t ^{\frac{1}{2}-\frac{d}{2p}}d\left( \frac{\tau}{t}\right)\cr
&\leq& 2C {\bf B}  \left(\frac{1}{2}-\frac{d}{2s},\frac{1}{2}- \frac{d}{2p}+\frac{d}{2s}\right).
\end{eqnarray}
For the case $t>1$, we can estimate
\begin{eqnarray}\label{Gp2}
	G_p(t)&\leq& 2C\int_0^{1/2} ( (t-\tau)^{-\frac{d}{2s}}+1) \tau^{-\frac{d}{2}
		(\frac{1}{p}-\frac{1}{s}+\frac{1}{d})}  e^{-\beta(t-\tau)} d\tau + 2C\int_{1/2}^1 ( (t-\tau)^{-\frac{d}{2s}}+1) \tau^{-\frac{d}{2}
		(\frac{1}{p}-\frac{1}{s}+\frac{1}{d})}  e^{-\beta(t-\tau)} d\tau\cr
&&+ 2C\int_{1}^t ( (t-\tau)^{-\frac{d}{2s}}+1)   e^{-\beta(t-\tau)} d\tau\cr
&\leq& 2C\int_0^{1/2} (2^{\frac{d}{2s}}+1) \tau^{-\frac{d}{2}
	(\frac{1}{p}-\frac{1}{s}+\frac{1}{d})}  d\tau + 2C\int_{1/2}^1 ( (t-\tau)^{-\frac{d}{2s}}+1) 2^{\frac{d}{2}
		(\frac{1}{p}-\frac{1}{s}+\frac{1}{d})}  e^{-\beta(t-\tau)} d\tau\cr
&&+ 2C\int_{1}^t ( (t-\tau)^{-\frac{d}{2s}}+1)    e^{-\beta(t-\tau)} d\tau\cr
&\leq& 2C (2^{\frac{d}{2s}}+1) \dfrac{2^{-\frac{1}{2}+ \frac{d}{2p}-\frac{d}{2s}}}{\frac{1}{2}- \frac{d}{2p}+\frac{d}{2s}} 
+ 2C\left(  2^{\frac{d}{2}
	(\frac{1}{p}-\frac{1}{s}+\frac{1}{d})}  +1\right) \int_{0}^t ( (t-\tau)^{-\frac{d}{2s}}+1)    e^{-\beta(t-\tau)} d\tau
\cr
	&\leq& 2C (2^{\frac{d}{2s}}+1) \dfrac{2^{-\frac{1}{2}+ \frac{d}{2p}-\frac{d}{2s}}}{\frac{1}{2}- \frac{d}{2p}+\frac{d}{2s}} 
+ 2C\left(  2^{\frac{d}{2}
	(\frac{1}{p}-\frac{1}{s}+\frac{1}{d})}  +1\right) \left[ \beta^{\frac{d}{2s}-1}{\bf \Gamma}  \left(1-\frac{d}{2s} \right)+\dfrac{1}{\beta}\right].
\end{eqnarray}
From inequalities \eqref{Gp1} and \eqref{Gp2}, we have 
$G_p(t)\leq K_p <+\infty.$

We remain to prove that
$$G_{s}(t):=2C \int_0^t [h_d(t)]^{-(\frac{1}{p}-\frac{1}{s}+\frac{1}{d})} [h_d(t-\tau)]^{\frac{1}{r}-\frac{1}{s}+\frac{1}{d}} [h_d(\tau)]^{\frac{1}{p}-\frac{1}{s}+\frac{1}{d}}  e^{-\beta(t-\tau)} d\tau\leq K_s<+\infty.$$
Indeed, for the case $0<t\leq 1$, we have
\begin{eqnarray}\label{Gs1}
	G_s(t)&\leq& 2C\int_0^t t^{\frac{d}{2}(\frac{1}{p}-\frac{1}{s}+\frac{1}{d})} (t-\tau)^{-\frac{d}{2}(\frac{1}{r}-\frac{1}{s}+\frac{1}{d})} \tau^{-\frac{d}{2}
		(\frac{1}{p}-\frac{1}{s}+\frac{1}{d})}  d\tau
	\cr&\leq& 2C\int_0^t  \left(1 -\frac{\tau}{t}\right) ^{-\frac{d}{2}(\frac{1}{p}+\frac{1}{d})} \left( \frac{\tau}{t}\right) ^{-\frac{d}{2}
		(\frac{1}{p}-\frac{1}{s}+\frac{1}{d})}  t ^{\frac{1}{2}-\frac{d}{2p}}d\left( \frac{\tau}{t}\right)
	\cr&\leq& 2C {\bf B}  \left(\frac{1}{2}-\frac{d}{2p},\frac{1}{2}- \frac{d}{2p}+\frac{d}{2s}\right).
\end{eqnarray}
For the case $t>1$, one has
\begin{eqnarray}\label{Gs2}
	G_s(t)&\leq& 2C\int_0^{1/2} ( (t-\tau)^{-\frac{d}{2}(\frac{1}{p}+\frac{1}{d})}+1) \tau^{-\frac{d}{2}
		(\frac{1}{p}-\frac{1}{s}+\frac{1}{d})}  e^{-\beta(t-\tau)} d\tau\cr
&& + 2C\int_{1/2}^1 ( (t-\tau)^{-\frac{d}{2}(\frac{1}{p}+\frac{1}{d})}+1) \tau^{-\frac{d}{2}
		(\frac{1}{p}-\frac{1}{s}+\frac{1}{d})}  e^{-\beta(t-\tau)} d\tau\cr
&&+  2C\int_{1}^t ( (t-\tau)^{-\frac{d}{2}(\frac{1}{p}+\frac{1}{d})}+1)   e^{-\beta(t-\tau)} d\tau\cr
&\leq& 2C\int_0^{1/2} (2^{\frac{d}{2p}+\frac{1}{2}}+1) \tau^{-\frac{d}{2}
		(\frac{1}{p}-\frac{1}{s}+\frac{1}{d})}  d\tau \cr
&&+ 2C\int_{1/2}^1 ( (t-\tau)^{-\frac{d}{2}(\frac{1}{p}+\frac{1}{d})}+1) 2^{\frac{d}{2}
		(\frac{1}{p}-\frac{1}{s}+\frac{1}{d})}  e^{-\beta(t-\tau)} d\tau\cr
&&+2C\int_{1}^t ( (t-\tau)^{-\frac{d}{2}(\frac{1}{p}+\frac{1}{d})}+1)    e^{-\beta(t-\tau)} d\tau\cr
&\leq&  2C (2^{\frac{d}{2p}+\frac{1}{2}}+1) \dfrac{2^{-\frac{1}{2}+ \frac{d}{2p}-\frac{d}{2s}}}{\frac{1}{2}- \frac{d}{2p}+\frac{d}{2s}} \cr
&&+2C\left(  2^{\frac{d}{2}
		(\frac{1}{p}-\frac{1}{s}+\frac{1}{d})}  +1\right) \int_{0}^t ( (t-\tau)^{-\frac{d}{2p}-\frac{1}{2}}+1)    e^{-\beta(t-\tau)} d\tau\cr
&\leq& 2C (2^{\frac{d}{2p}+\frac{1}{2}}+1) \dfrac{2^{-\frac{1}{2}+ \frac{d}{2p}-\frac{d}{2s}}}{\frac{1}{2}- \frac{d}{2p}+\frac{d}{2s}} 
	+ 2C\left(  2^{\frac{d}{2}
		(\frac{1}{p}-\frac{1}{s}+\frac{1}{d})}  +1\right) \left[ \beta^{\frac{d}{2p}-\frac{1}{2}}{\bf \Gamma}  \left(\frac{1}{2}-\frac{d}{2p} \right)+\dfrac{1}{\beta}\right] .
\end{eqnarray}
Combining inequalities \eqref{Gs1} and \eqref{Gs2}, we get
$G_s(t)\leq K_s <+\infty.$

\noindent
{\bf Conflict of interest:} The authors declare that there is no conflict of interest\\
{\bf Data availability statement:} There is no associated data and no new data was created in this study


\begin{thebibliography}{9}
\bibitem{An2009} J.-P. Anker and V. Pierfelice, {\it Nonlinear Schr\"odinger equation on real hyperbolic spaces}, Ann. I. H. Poincaré – AN 26 (2009), pages 1853–1869

\bibitem{Al2011} M.F.D. Almeida, L.C.F. Ferreira, {\it On the well-posedness and large time behavior for Boussineq equations in Morrey spaces}, Differential integral equations 24 (7-8) (2011) 667-684.

\bibitem{Ba} L. Brandolese and M.E. Schonbek, {\it Large time decay and growth for solutions of a viscous Boussineq system}, Trans. Amer. Math. Soc. 364 (10) (2012) 5057-5090.



\bibitem{Barrow2020} J.D. Barrow, {\it Non-Euclidean Newtonian cosmology}, Classical and Quantum Gravity, Vol. 37, Num. 12 (2020)

\bibitem {Br2012}L. Brandolese and M.E. Schonbek, \textit{Large time decay and
growth for solutions of a viscous Boussinesq system}, Trans. Amer. Math. Soc.
\textbf{364} (10) (2012), 5057-5090.

\bibitem {Br2020}L. Brandolese and J. He, \textit{Uniqueness theorems for the
Boussinesq system}, Tohoku Math. J. 72 (2) (2020), 283-297.

\bibitem{Ca1980} J.R. Cannon and E. DiBenedetto, {\it The initial value problem for the Boussinesq equations with data in $L^p$, in Approximation Methods for Navier–Stokes Problems}, Edited by Rautmann, R., Lect. Notes in Math.,
Springer-Verlag, Berlin, 771, 1980.

\bibitem{Cao1999} C. Cao, M. Rammaha, and E. Titi, {\it The Navier-Stokes equations on the rotating $2-$D sphere: Gevrey regularity and asymptotic degrees of freedom}, In: Z. Angew. Math. Phys. 50.3 (1999), pp. 341-360.

\bibitem {Chandra}S. Chandrasekhar, \textit{Hydrodynamic and Hydromagnetic Stability}, Dover, New York, 1981.

\bibitem{Cz1} M. Czubak and C.H. Chan, {\it Non-uniqueness of the Leray–Hopf solutions in the hyperbolic setting}, Dyn. PDE 10(1), 43–77 (2013).

\bibitem{Cz2} M. Czubak and C.H. Chan, {\it Remarks on the weak formulation of the Navier–Stokes equations on the $2$D-hyperbolic space}, Annales de l'Institut Henri Poincare (C) Non Linear Analysis, Volume 33, Issue 3, May-June 2016, Pages 655-698.

\bibitem{Cz3} M. Czubak, C.H. Chan and M. Disconzi, {\it The formulation of the Navier-Stokes equations on Riemannian manifolds} In: J. Geom. Phys. 121 (2017), pp. 335–346.

\bibitem{DalKre} Ju. L. Daleckii, M. G. Krein, {\it Stability of Solutions of Differential Equations in Banach Spaces}, Transl. Amer. Math. Soc.  Provindence RI, {1974.}

\bibitem{Da} E. Damek and F. Ricci, {\it A class of nonsymmetric harmonic Riemannian spaces}, Bull. Amer. Math. Soc. 27 (1992), 139-142.

\bibitem {Danchin2009}R. Danchin and M. Paicu, \textit{Global well-posedness
issue for the inviscid Boussinesq system with Youdovitch's type data}, Commun.
Math. Phys. 290 (2009), 1-14.

\bibitem {Danchin2008}R. Danchin and M. Paicu, \textit{Existence and
uniqueness results for the Boussinesq system with data in Lorentz spaces}, Phys. D 237 (10-12) (2008), 1444-1460.

\bibitem{Do} H. Dong and Q. S. Zhang, {\it Time analyticity for the heat equation and Navier-Stokes equations}, Journal of Functional Analysis 279(4):108563, April 2020, DOI: 10.1016/j.jfa.2020.108563

\bibitem{EbiMa} D.G. Ebin, J.E. Marsden, {\it Groups of diffeomorphisms and the motion of an incompressible fluid}, Ann. of Math. (2) 92 (1970), 102-163.

\bibitem{Er} P. Erbelein, {\it Geometry of non positively curved manifolds}, Chicago Lectures in Mathematics, 449, (1996).

\bibitem{Fa2018} S. Fang and D. Luo, {\it Constantin and Iyer’s representation formula for the Navier-Stokes equations on manifolds}, Potential Analysis, 48: 181-206 (2018).

\bibitem{Fa2020} S. Fang, {\it Nash Embedding, Shape Operator and Navier-Stokes Equation on a Riemannian Manifold}, Acta Mathematicae Applicatae Sinica, English Series Vol. 36, No. 2 (2020) 237-252.

\bibitem{Fe2006} L.C.F. Ferreira and E.J. Villamizar-Roa, {\it Well-posedness and asymptotic behaviour for the convection
problem in $\mathbb{R}^n$}, Nonlinearity, 19, 2169-2191 (2006).

\bibitem{Fe2008} L.C.F. Ferreira and E.J. Villamizar-Roa, {\it Existence of solutions to the convection problem in a pseudomeasure-type space}, Proc. R. Soc. Lond. Ser. A Math. Phys. Eng. Sci. 464, no. 2096,
1983-1999, (2008).

\bibitem{Fe2010} L.C.F. Ferreira and E.J. Villamizar-Roa, {\it On the stability problem for the Boussinesq equations in weak-$L^p$spaces}, Commun. Pure Appl. Anal. 9 (3) (2010) 667–684.

\bibitem {Fe2016} L.C.F. Ferreira, \textit{On a bilinear estimate in weak-Morrey spaces and uniqueness for Navier-Stokes equations}, J. Math. Pures Appl. 105 (2) (2016), 228-247.

\bibitem{FePha2023} L.C.F. Ferreira and P.T. Xuan, {\it On uniqueness of mild solutions for Boussinesq equations in Morrey-type spaces}, Applied Mathematics Letters, Vol. 137, Num. 2 (2023) 1-8.

\bibitem {Fi1969}P.C. Fife and D.D. Joseph, \textit{Existence of convective
solutions of the generalized Bernard problem which are analytic in their
norm}, Arch. Rational Mech. Anal. 33 (1969), 116-138.

\bibitem{FuKa} H. Fujita and T. Kato, {\it On the Navier-Stokes initial value problem}, I, Arch. Rat. Mech. Anal. 16,(1961), 269-315.

\bibitem{Go} A.A. Grigor\'yan, {\it Analytic and geometric background of recurrence and non-explostion of the Brownian motion on Riemannian manifolds}, Bulletin (New Series) of the American mathematical society,
Vol. 36, Num. 2, Pages 135–249 (1999)

\bibitem{He} E. Hebey, {\it Nonlinear Analysis on Manifolds: Sobolev Spaces and Inequalities, Courant Lectures in Mathematics}, New York (2000) AMS.

\bibitem{Hel} S. Helgason, {\it Geometric Analysis on Symmetric Spaces}, Amer. Math. Soc., 1994.

\bibitem {Hi1995}T. Hishida, \textit{Global Existence and Exponential
Stability of Convection}, J. Math. Anal. Appl. 196 (2) (1995), 699-721.

\bibitem {Hi1997}T. Hishida, \textit{On a class of Stable Steady flow to the
Exterior Convection Problem}, Journal of Differential Equations 141 (1)
(1997), 54-85.


\bibitem{HuyXuan2021} N.T. Huy, P.T. Xuan, V.T.N. Ha and V.T. Mai, {\it Periodic Solutions of Navier-Stokes Equations on Non-compact Einstein Manifolds with Negative Ricci Curvature Tensor}, Analysis and Mathematical Physics (2021) 11:60.

\bibitem{HuyXuan2022} N.T. Huy, P.T. Xuan, V.T.N. Ha, N.T. Van, {\it On periodic solutions for the incompressible Navier-Stokes equations on non-compact manifolds}, Taiwan. J. Math., 26 (3) (June 2022), pp. 607-633, 10.11650/tjm/211205

\bibitem{HuyXuan22} N.T. Huy, P.T. Xuan, V.T.N. Ha and L.T. Sac, {\it Existence and Stability of Periodic and Almost Periodic Solutions to the Boussinesq System in Unbounded Domains}, Acta. Math. Sci, Vol. 42, Iss. 5, 1875–1901 (2022).


\bibitem{HuyVan2022} N.T. Huy, V.T.N. Ha and N.T. Van, {\it Stability and periodicity of solutions to Navier-Stokes equations on non-compact riemannian manifolds with negative curvature}, Anal.Math.Phys. 12, {\bf 89} (2022).

\bibitem{Khe} B. Khesin and G. Misiolek, {\it The Euler and Navier-Stokes equations on the hyperbolic plane}, Proc. Natl. Acad. Sci. 109 (2012), pages 18324-18326

\bibitem{Il1991} A.A. Il'yin, {\it The Navier-Stokes and Euler equations on two-dimensional closed manifolds}, Mat. Sb. 181 (1990), 521-539; English transl. in Math. USSR Sb. 69 (1991).


\bibitem{Jo} J. Jost, {\it Riemannian geometry and geometric analysis}, Universitext. Springer-Verlag, Berlin, fifth edition (2008).

\bibitem{Ko2008} M. Kobayashi, {\it On the Navier-Stokes equations on manifolds with curvature}, In: J. Engrg. Math. 60.1 (2008), pp. 55-68.

\bibitem {Komo2015}C. Komo, \textit{Uniqueness criteria and strong solutions of the Boussinesq equations in completely general domains}, Z. Anal. Anwend. 34 (2) (2015), 147-164.

\bibitem{Li2016} L. A. Lichtenfelz, {\it Nonuniqueness of solutions of the Navier–Stokes equations on Riemannian manifolds}, Annals of Global Analysis and Geometry volume 50, pages237–248(2016).

\bibitem {Li-Wang2021}Z. Li and W. Wang, \textit{Norm inflation for the Boussinesq system}, Discrete Contin. Dyn. Syst. Ser. B 26 (10) (2021), 5449-5463.

\bibitem {Liu2014}X. Liu and Y. Li, \textit{On the stability of global
solutions to the 3D Boussinesq system}, Nonlinear Analysis 95 (2014), 580-591.

\bibitem{Loho} N. Lohou\'e, {\it Estimation des projecteurs de De Rham Hodge de certaines vari\'et\'e riemanniennes non compactes}, Math. Nachr. 279, 3 (2006), 272-298.


\bibitem{MiTa2001} M. Mitrea and M. Taylor, {\it Navier-Stokes equations on Lipschitz domains in Riemannian manifolds}, Math Ann 321, 955-987 (2001). 

\bibitem {Mo1991}H. Morimoto, \textit{Non-stationary Boussinesq equations},
Proc. Japan Acad. Ser. A math. Sci. 67 (5) (1991), 159-161.

\bibitem {Na2020}K. Nakao, \textit{On time-periodic solutions to the
Boussinesq equations in exterior domains}, J. Math. Anal. Appl. 482 (2)
(2020), 123537, 16 pages.


\bibitem{Pa} J. Pauqert, {\it Introduction to hyperboloid geometry}, lecture note, (2016).

\bibitem{Pi} V. Pierfelice, {\it The incompressible Navier-Stokes equations on non-compact manifolds}, Journal of Geometric Analysis, 27(1) (2017), 577-617.

\bibitem{Re} W.F. Reynords, {\it Hyperbolic Geometry on a Hyperboloid}, The American Mathematical Monthly, Vol. 100, No. 5 (May, 1993), pp. 442-455.

\bibitem{Sa} M. Samavaki and J. Tuomela, {\it Navier–Stokes equations on Riemannian manifolds}, Journal of Geometry and Physics 148 (2020) 103543

\bibitem{Tay}  M. Taylor, {\it Partial Differential Equations III: Nonlinear equations},volume 117 of Applied Mathematical Sciences. Springer New York second edition (2011).

\bibitem{Te01} R. Temam, {\it Navier-Stokes equations}, AMS Chelsea Publishing, Providence, RI, 2001.

\bibitem{XVQ2023} P.T. Xuan, N.T. Van and B. Quoc, {\it On Asymptotically Almost Periodic Solution of Parabolic Equations on real hyperbolic Manifolds}, Journal of Mathematical Analysis and Applications, Vol. 517, Iss. 1 (2023), 126578.

\bibitem{XuanVan2021} P.T. Xuan and N.T. Van, {\it On asymptotically almost periodic solutions to the Navier-Stokes equations on hyperbolic manifolds}, Journal of Fixed Point Theory and Applications, Vol. 25, No. 71 (2023), 33 pages.

\bibitem{XVT2024} P.T. Xuan, N.T. Van and T.V. Thuy, {\it Periodic solutions for Boussinesq systems in weak-Morrey spaces}, Journal of Mathematical Analysis and Applications, Vol. 537, No. 1 (2024), 128255.

\bibitem{XuanTrung2022} P.T. Xuan, T.T. Ngoc N.T.Van and  H. Trung, {\it On periodic solution for the Boussinesq system on real hyperbolic manifolds}, Dynamical Systems: An International Journal (2024), 26 pages. https://doi.org/10.1080/14689367.2024.2427791

\bibitem{NgocXuan2024} P.T. Xuan and T.T. Ngoc, {\it On pseudo almost periodic solutions for Boussinesq systems on real hyperbolic Manifolds}, accepted for publication in Filomat Journal (2024), 14 pages.

\bibitem{XN2024} P.T. Xuan and T.T. Ngoc, {\it Stability of Solutions of Stationary Boussinesq Systems in Weak-Morrey Spaces}, Bull Braz Math Soc, New Series {\bf 55}, 49 (2024).



\bibitem{Zha} Q. S. Zhang, {\it The ill-posed Navier-Stokes equation on connected sums of $\mathbb{R}^3$}, Complex Variables and Elliptic Equations, 51:8-11, 1059-1063 (2006).

\end{thebibliography}
\end{document}